\documentclass[notitlepage,11pt,draft]{article}
\usepackage{amssymb,amsfonts,amsmath,geometry,esint,comment,upgreek,cases,graphicx,bbm,ulem}
\usepackage{scalerel,stackengine,mathtools,stmaryrd,soul,array}

\usepackage{color}

\catcode`\@=11
\@addtoreset{equation}{section}

\catcode`\@=12

\allowdisplaybreaks

\def\eps{\varepsilon}

\def\R{\mathbb{R}}

\def\f{\varphi}

\def\S{\mathbb S}
\def\rmH{{\rm H}}

\def\ird{\int\limits_{\R^d}}

\def\tgamma{\widetilde{\gamma}}

\def\pstar{{p^*}}
\def\twostar{{2^*}}

\def\Ps{|\Pi\sigma|}

\def\div{{\rm div}}

\def\proof{\noindent{\textbf{Proof. }}}
\def\QED{\hfill {$\square$}\goodbreak \medskip}

\newtheorem{Theorem}{Theorem}[section]
\newtheorem{Lemma}[Theorem]{Lemma}
\newtheorem{Proposition}[Theorem]{Proposition}
\newtheorem{Corollary}[Theorem]{Corollary}
\newtheorem{Remark}[Theorem]{Remark}

\begin{document}

\title{
Hardy-Sobolev inequalities involving \\
mixed radially and cylindrically symmetric weights}

\author{
Gabriele Cora\footnote{Département de Mathématique
Université Libre de Bruxelles, Belgium. Email: {gabriele.cora@ulb.be}, orcid.org/0000-0002-0090-5470.},~~
Roberta Musina\footnote{DMIF, Universit\`a di Udine, Italy.
Email: {roberta.musina@uniud.it}, orcid.org/0000-0003-4835-8004.},~~
Alexander I. Nazarov\footnote{
St. Petersburg Department of Steklov Institute, St{.} Petersburg, Russia. E-mail: al.il.nazarov@gmail.com, 
orcid.org/0000-0001-9174-7000.}
}

\date{}

\maketitle

\centerline{\it Dedicated to Prof. Gianni Mancini}

\bigskip
\noindent
{\small {\bf Abstract.} 
We deal with weighted Hardy--Sobolev type inequalities for functions on $\R^d$, $d\ge 2$. The weights involved are anisotropic, given by products of powers of the distance to the origin and to a  nontrivial subspace.  
We establish necessary and sufficient conditions for validity of these inequalities, and investigate the existence/nonexistence of extremal functions.

\medskip



\section{Introduction}

We take integers $d>k\ge 1$ and identify
$
\R^d\equiv \R^{d-k}\times \R^k\ni (x,y)\equiv z~\!.
$
Given $p>1$ and $a,b\in\R$, we study dilation-invariant inequalities related to the differential operators
\begin{equation}
\label{eq:operator}
\mathcal L_{a,b}(u):=-\div(|y|^a|z|^{-b}|\nabla u|^{p-2}\nabla u)~\!,
\end{equation}
which
arise in a variety of scenarios. Below are a few motivating examples.

\medskip

\noindent$-$ When $k=1$ and $a\in(-1,1)$, the elliptic operator
$-\div(|y|^{a}\nabla u)$
is related to the fractional Laplacian of order $s=\frac{1-a}{2}$ on $\R^{d-1}$, see the milestone paper \cite{CS} by Caffarelli and Silvestre.

\medskip

\noindent$-$ 
Hardy and Hardy--Sobolev type inequalities play a fundamental role in establishing regularity results for variational problems driven by operators of the form
$-\div(|y|^aA(z)\nabla u)$.  See, for instance,~\cite{CV, DFM, DM, Fio, STV1, STV2} and references therein.

\medskip

\noindent$-$
For $k = d-1$, 
the asymptotic stability of solutions to the Navier--Stokes equations has been 
investigated in  \cite{YYLi1, YYLi}  by employing anisotropic integral inequalities related to the functional
$$
\langle \mathcal L_{a,b}(u),u\rangle = \int\limits_{\R^d}|y|^{a}|z|^{-b}|\nabla u|^p~\!dz~\!.
$$

To formulate the problems under consideration, we introduce, for any
$t\in\R$, the quantity
$$
\rmH_{t}=\frac{d-p+t}{p}~\!.
$$
Assume  
\begin{equation}
\label{eq:CMN_assu}
k+a>0~,\qquad b<d-p+a=p\rmH_a~\!,
\end{equation}
so that $\rmH_{a-b}>0$. 
It is easy to show, see for instance \cite[Theorem 1]{CMN}, that the Hardy type inequality
\begin{equation}
\label{eq:CMN}
\rmH_{a-b}^p\ird |y|^a|z|^{-b-p}|u|^p~\!dz\le \ird |y|^a|z|^{-b}|\nabla u|^p~\!dz
\end{equation}
holds with a sharp constant $\rmH_{a-b}^p$ for any $u\in \mathcal C^\infty_c(\R^d)$. 
We emphasize that the conditions in (\ref{eq:CMN_assu}) are needed to have the local integrability of the weights 
involved.

Thanks to (\ref{eq:CMN}), we can define
the reflexive Banach space $$\mathcal D^{1,p}(\R^d;|y|^{a}|z|^{-b}~\!dz),$$
as the completion of ${\cal C}^\infty_c(\R^d)$ 
with respect to the norm $\||y|^{\frac{a}{p}}|z|^{-\frac{b}{p}}\nabla \cdot\|_{L^p}$. Then 
(\ref{eq:CMN}) holds for  any  $u\in \mathcal D^{1,p}(\R^d;|y|^a|z|^{-b}dz)$.
More general Hardy type inequalities for functions in $\mathcal D^{1,p}(\R^d;|y|^{a}|z|^{-b}~\!dz)$
were proved in \cite{MN_H},  see also~\cite{YYLi1, CMN, HY}. 

\medskip

Throughout this paper we take $q>p$. We seek  conditions on $q$ and  on the new exponents 
$\eta, \gamma\in \R$ to ensure that the inequality
$$
c\Big(\displaystyle\int_{\R^d}|y|^{\eta\frac{q}{p}}|z|^{-\gamma\frac{q}{p}}|u|^q~\!dz\Big)^\frac{p}{q}
\le \int_{\R^d}|y|^{a}|z|^{-b}|\nabla u|^p~\!dz~,\qquad u\in \mathcal D^{1,p}(\R^d;|y|^{a}|z|^{-b} dz),
$$
holds for some  constant $c>0$. 
By taking into account the
action of dilations we are forced to take $\eta\frac{q}{p}=-d+q\rmH_{a-b+\gamma}$. Further, to have the local integrability  of the weight in the left hand side,
we need to impose
\begin{equation}
\label{eq:CMN_assu2}
q\rmH_{a-b+\gamma}>d-k~\!. 
\end{equation}

We are therefore naturally lead to consider the optimal constant
\begin{equation}
\label{eq:HS_constant}
S_{a,b,\gamma}(q)= \inf_{u\in \mathcal D^{1,p}(\R^d;|y|^{a}|z|^{-b} dz)}
\frac{\displaystyle\int_{\R^d}|y|^{a}|z|^{-b}|\nabla u|^p~\!dz}
{\Big(\displaystyle\int_{\R^d}|y|^{-d+q\rmH_{a-b+\gamma}}|z|^{-\gamma\frac{q}{p}}|u|^q~\!dz\Big)^\frac{p}{q}}~\!,
\end{equation}
under the  additional assumption (\ref{eq:CMN_assu2}).

The infimum $S_{a,b,\gamma}(q)$ has been extensively studied for several specific choices of the parameters:

\medskip

\noindent$-$ 
\textit{Purely spherical weights.} This case is related to the classical Il'in--Caffarelli--Kohn--Nirenberg inequalities~\cite{il, CKN}, and corresponds to the choice $a=0$ and $q\rmH_{a-b+\gamma}=d$. 
The existence of extremals was first studied in \cite{CW} in the Hilbertian case $p=2$ and then in \cite{GM} for $p\in(1,\infty)$.

If $p<d$, by further setting $\gamma=b=0$ and $\pstar=\frac{dp}{d-p}$, we recover the (unweighted) Sobolev constant
$$
S= \inf_{u\in \mathcal D^{1,p}(\R^d)} \frac{\displaystyle\int_{\R^d}|\nabla u|^p~\!dz}
{\Big(\displaystyle\int_{\R^d}|u|^{p^*}~\!dz\Big)^\frac{p}{p^*}}~\!,
$$
whose explicit value was computed by Aubin \cite{Au} and Talenti \cite{Ta}.

\medskip

\noindent$-$ \textit{Purely cylindrical weights.}  If $\gamma=b=0$, then the  infimum $S_{a,0,0}(q)$ reduces to the {\it Maz'ya constant} 
\begin{equation}
\label{eq:Ma}
M_{a}(q)
= \inf_{u\in \mathcal D^{1,p}(\R^d;|y|^{a} dz)}
\frac{\displaystyle\int_{\R^d}|y|^{a}|\nabla u|^p~\!dz}
{\Big(\displaystyle\int_{\R^d}|y|^{-d+q\rmH_a}|u|^q~\!dz\Big)^\frac{p}{q}}~\!,
\end{equation}
see \cite[Corollary 2.1.7/1]{Ma}. Together with the Sobolev constant, $M_a(q)$ plays a crucial role in the present paper.
It turns out that $M_a(q)>0$ if and only if 
\begin{equation}
\label{eq:assu_cyl}
k+a>0~,
\qquad q\rmH_a>d-k~, \qquad (d-p)q\le dp~\!.
\end{equation}
Notice that the first two  inequalities in (\ref{eq:assu_cyl}) are equivalent to (\ref{eq:CMN_assu}), (\ref{eq:CMN_assu2}) when $\gamma=b=0$;
the last one takes into account the Sobolev embedding theorem.

The problem of the existence of extremals for $M_{a}(q)$ was first addressed in \cite{BT} in case
$a=0$, $p<k$ and $q<\pstar$. It was then studied in \cite{TT, Mu} for $p=2$, and
later treated in generality in \cite{GM}.

\medskip

\noindent$-$ 
\textit{The case $k=d-1$.}
This case is explored in \cite{YYLi1}  
and more specifically in \cite{YYLi}, within the context of more general integral inequalities.

\bigskip

Our first result  follows.

\begin{Theorem}
\label{T:Main1}
Let  
(\ref{eq:CMN_assu}), (\ref{eq:CMN_assu2}) hold. Then $S_{a,b,\gamma}(q)>0$ if and only if the following conditions are satisfied:
\begin{equation}
\label{eq:assu}
a1)~ (d-p)q\le dp~,\qquad a2)~\gamma\ge b.
\end{equation}
If in addition $p<d$ then
\begin{equation}
\label{eq:S_critical}
S_{a,b,\gamma}(p^*)= \inf_{u\in \mathcal D^{1,p}(\R^d;|y|^{a}|z|^{-b} dz)}
\frac{\displaystyle\int_{\R^d}|y|^{a}|z|^{-b}|\nabla u|^p~\!dz}
{\Big(\displaystyle\int_{\R^d}|y|^{(a-b+\gamma)\frac{\pstar}{p}}|z|^{-\gamma\frac{\pstar}{p}}|u|^\pstar~\!dz\Big)^\frac{p}{\pstar}}\le S.
\end{equation}
\end{Theorem}

\medskip

It is evident that the assumption $a1)$ holds whenever $p\ge d$. If $p< d$, then the restriction 
$q\le p^*$ is necessary in order to ensure that $S_{a,b,\gamma}(q)>0$, due to the Sobolev embedding theorem.

\medskip

Assumption $a2)$ prevents vanishing along the {\it singular set} 
\begin{equation}
\label{eq:Sigma0}
\Sigma_0=\{y=0\}~\!,
\end{equation}
a phenomenon already observed in \cite{YYLi} (for $k=d-1$) and in \cite{MN_H}.

It is worth noting that, given $a,b$ as in (\ref{eq:CMN_assu}) and $q>p$ satisfying 
$a1)$, the function 
$$
\gamma\mapsto S_{a,b,\gamma}(q)~,\qquad \gamma\in \big(b+\tfrac{p}{q}(d-k-q\rmH_{a}),\infty\big)\cap[b,\infty)
$$
is (positive and)  non decreasing.

Notice that the {\it bottom case} $\gamma=b$ can only occur if $q\rmH_a>d-k$, which requires $p\rmH_a=d-p+a>0$.
In this case, it trivially holds that $S_{a,b,\gamma}(q)\ge S_{a,b,b}(q)$ for any $\gamma\ge b$.

\medskip

\medskip

We next address  the attainability of the optimal constant $S_{a,b,\gamma}(q)$.

The minimization problem in \eqref{eq:HS_constant} is clearly non-compact.
By the invariance of the Rayleigh quotient in \eqref{eq:HS_constant} under rescaling,
one can easily construct minimizing sequences concentrating either at the origin
or at infinity.
More subtle and more severe loss of compactness phenomena occur in the critical case
$q=\pstar$ (due to the action of the translation group in $\R^d$),
and in the bottom case $\gamma=b$ (due to translations along the singular set $\Sigma_0$).

The next table summarizes our existence results\footnote{We agree that (\ref{eq:CMN_assu}) and (\ref{eq:CMN_assu2}) hold,
and that $\pstar=\infty$ if $p\ge d$.}
$$
\begin{tabular}{|c|c|c|}
\hline
{\it Assumptions}&{\it Lack of compactness due to}&{\it $S_{a,b,\gamma}(q)$ is achieved \dots}\\
\hline\hline
$p<q<\pstar$ & dilations in $\R^d$&\dots {\it always}\\ 
and $\gamma>b$&&see $i)$ in Theorem \ref{T:Main2} \\ 
\hline 
$q=\pstar$& dilations and  translations in $\R^d$&\dots {\it if $S_{a,b,\gamma}(\pstar)<S$}\\
and $\gamma>b$&&see $ii)$ in Theorem \ref{T:Main2} \\
\hline
$p<q\le \pstar$& dilations in $\R^d$ &\dots {\it if $S_{a,b,b}(q)<M_a(q)$}\\
and $\gamma=b$&and translations along $\Sigma_0=\{y=0\}$&see Theorem \ref{T:Main3}\\ 
\hline
\end{tabular}
$$
We emphasize that nonexistence may indeed occur in the critical case when
$S_{a,b,\gamma}(\pstar)=S$, as well as in the bottom case when
$S_{a,b,b}(q)=M_a(q)$.
This phenomenon has already been observed in both the purely spherical
and the purely cylindrical settings (see also Theorem~\ref{T:H_bottom_sub}).

We are in position to state our existence theorems in a  rigorous way.

\begin{Theorem} 
\label{T:Main2}
Let (\ref{eq:CMN_assu}) and (\ref{eq:CMN_assu2}) be satisfied. In addition, assume $\gamma > b$. Then
\begin{itemize}
\item[$i)$] If $p\ge d$ or $q<p^*$, then $S_{a,b,\gamma}(q)$ is achieved;
\item[$ii)$] If $p<d$ and $S_{a,b,\gamma}(p^*)< S$, then $S_{a,b,\gamma}(p^*)$ is achieved.
\end{itemize}
\end{Theorem}

\medskip 
We outline the main steps in the proof of Theorem \ref{T:Main2}.
To overcome the lack of compactness caused by dilations in $\R^d$, we  use 
Ekeland’s variational principle together with  an appropriate concentration function
to  construct a well-behaved minimizing sequence
for problem (\ref{eq:HS_constant}). In the critical case $q=\pstar$, the additional loss of compactness due to translations is 
handled through a Sacks--Uhlenbeck type strategy, inspired by \cite{SU}
and based on an {\it $\eps$-compactness criterion}. 

\medskip

Next we deal with the bottom case $\gamma=b$.  
Assuming that $a, q$ satisfy (\ref{eq:assu_cyl}), we regard at $S_{a,b,b}(q)$
as a function of $b\in (-\infty,p\rmH_a)$, compare with (\ref{eq:CMN_assu}). Explicitly, we have 
\begin{equation}
\label{eq:Mb}
S_{a,b,b}(q)=\inf_{u\in \mathcal D^{1,p}(\R^d;|y|^{a}|z|^{-b} dz)}
\frac{\displaystyle\int_{\R^d}|y|^{a}|z|^{-b}|\nabla u|^p~\!dz}
{\Big(\displaystyle\int_{\R^d}|y|^{-d+q\rmH_a}|z|^{-b\frac{q}{p}}|u|^q~\!dz\Big)^\frac{p}{q}}~\!.
\end{equation}
Notice that for $b=0$, the infimum  $S_{a,0,0}(q)$  reduces to the Maz'ya constant in (\ref{eq:Ma}).

\begin{Theorem} [Bottom case]
\label{T:Main3}
Let (\ref{eq:assu_cyl}) be satisfied. For every $b<p\rmH_a$ the following facts hold:
\begin{itemize}
\item[$i)$] $S_{a,b,b}(q)\le M_{a}(q)$;
\item[$ii)$] If $S_{a,b,b}(q)< M_{a}(q)$, then $S_{a,b,b}(q)$ is achieved.
\end{itemize}
\end{Theorem}

\medskip
The proof of Theorem \ref{T:Main3} follows a Sacks--Uhlenbeck type strategy similar to that of $ii)$ in Theorem \ref{T:Main2}, but it relies on a different concentration 
function and needs further energy estimates on the ratio in (\ref{eq:Mb}) to prevent vanishing phenomena.

\medskip

For $p=2$ we further analyze the minimization problem (\ref{eq:HS_constant}) in the bottom case $\gamma=b$. Namely, we deal with
\begin{equation}
\label{eq:Mb2}
S_{a,b,b}(q)=\inf_{u\in \mathcal D^{1,2}(\R^d;|y|^{a}|z|^{-b} dz)}
\frac{\displaystyle\int_{\R^d}|y|^{a}|z|^{-b}|\nabla u|^2~\!dz}
{\Big(\displaystyle\int_{\R^d}|y|^{-d+q\rmH_a}|z|^{-b\frac{q}{2}}|u|^q~\!dz\Big)^\frac{2}{q}}~\!.
\end{equation}
The underlying Hilbertian structure allows us to use a natural functional transform, see (\ref{eq:Ta}), to  rewrite $S_{a,b,b}(q)$ in the form
$$
\inf_{v\in \mathcal D^{1,2}(\R^d;|y|^{a}dz) }
\frac{\displaystyle~
\int_{\R^d}|y|^{a}|\nabla v|^2~\!dz+\tfrac{b}{2}\big(\tfrac{b}{2}-2\rmH_a\big)\int_{\R^d}|y|^{a}|z|^{-2}|v|^2~\!dz~}
{\displaystyle{\Big(
\int_{\R^d}|y|^{-d+q\rmH_{a}}|v|^q~\!dz\Big)^\frac{2}{q}}}~\!.
$$
In this way we  obtain a more complete picture of existence and nonexistence of extremals. 
Here we restrict ourselves to stating our results for subcritical exponents $q>2$. 
A more general theorem, which also covers the critical exponent $q=2^*$, is presented in Theorem \ref{T:H_bottom_sub} of Section \ref{S:Hilbert}.

\begin{Theorem}
\label{T:MainH}
Let $p=2<q<2^*$ (we agree that $2^*=\infty$ if $d=2$). Assume 
$$
a+k\frac{2}{2^*}>2(d-k)\Big(\frac{1}{q}-\frac{1}{2^*}\Big)~\!,
$$
compare with (\ref{eq:assu_cyl}).
Then the following facts hold:
\begin{itemize}
\item[$i)$]  
$\begin{cases}
S_{a,b,b}(q)= M_{a}(q)&\text{if $b<0$, and in this case $S_{a,b,b}(q)$  is not achieved;}\\
S_{a,b,b}(q)< M_{a}(q)&\text{if $0<b<2\rmH_a$,  and in this case $S_{a,b,b}(q)$  is achieved.}
\end{cases}
$

\smallskip
(It is known that the Maz'ya constant $M_a(q)=S_{a,0,0}(q)$ is  achieved);
\item[$ii)$] The map $b\mapsto S_{a,b,b}(q)$ is strictly decreasing and continuous on $[0,2\rmH_a)$.
\end{itemize}
\end{Theorem}

We point out that the  existence result in the recent paper \cite{2025}, where  $p=2<d$ and $q<2^*$ are assumed,
is incorrect.

\medskip

The paper is organized as follows. 
In the next section, we prove Theorem \ref{T:Main1}. The existence Theorems \ref{T:Main2} and \ref{T:Main3} are 
established in Section \ref{S:Main23}, 
after some preliminary results.

Section~\ref{S:strict} provides sufficient conditions to have 
$S_{a,b,\gamma}(\pstar)<S$ and
$S_{a,b,\gamma}(q)<M_a(q)$, along with some nonexistence results.

The Hilbertian, bottom  case ($p=2$ and $\gamma=b$) is addressed in Section \ref{S:Hilbert}. 

Appendix~\ref{A:Examples} contains few illustrative examples. 

Several questions concerning the minimization problem (\ref{eq:Ma}) remain open. The relevant known results 
used in the present paper are summarized in Appendix \ref{A}.

\bigskip

Through the paper, any positive constant whose value is not important is denoted 
by $c$. It may take different values at different places. To indicate that a 
constant depends on some parameters 
we list them in parentheses.

\section{Proof of Theorem \ref{T:Main1}}
\label{S:Main1}

As already noted in the introduction, assumption $a1)$ in (\ref{eq:assu}) is necessary to have $S_{a,b,\gamma}(q)>0$. We now show that $a2)$ 
 is also required.

Take $x_0\in\R^{d-k}\setminus\{0\}$ with $|x_0|=1$ and put $z_0=(x_0,0)$. Take a nontrivial function ${u}\in \mathcal C^\infty_c(\R^d)$ and for any 
integer $h> 1$ put
$$
{u}_h(z)={u}(z-hz_0).
$$
We have 
$$
\frac{\displaystyle \int_{\R^d}|y|^{a}|z|^{-b}|\nabla {u}_h|^p~\!dz}
{\Big(\displaystyle \int_{\R^d}|y|^{-d+q\rmH_{a-b+\gamma}}|z|^{-\gamma\frac{q}{p}}|{u}_h|^q~\!dz\Big)^{\frac{p}{q}}}
=h^{\gamma-b} ~
\frac{\displaystyle \int_{\R^d}|y|^{a}|h^{-1}z+z_0|^{-b}|\nabla {u}|^p~\!dz}
{\Big(\displaystyle \int_{\R^d}|y|^{-d+q\rmH_{a-b+\gamma}}|h^{-1}z+z_0|^{-\gamma\frac{q}{p}}|{u}|^q~\!dz\Big)^{\frac{p}{q}}}
$$
so that
\begin{equation}
\label{eq:prima}
\frac{\displaystyle \int_{\R^d}|y|^{a}|z|^{-b}|\nabla {u}_h|^p~\!dz}
{\Big(\displaystyle \int_{\R^d}|y|^{-d+q\rmH_{a-b+\gamma}}|z|^{-\gamma\frac{q}{p}}|{u}_h|^q~\!dz\Big)^{\frac{p}{q}}}
=h^{\gamma-b} 
\Bigg[\frac{\displaystyle \int_{\R^d}|y|^{a}|\nabla {u}|^p~\!dz}
{\Big(\displaystyle \int_{\R^d}|y|^{-d+q\rmH_{a-b+\gamma}}|{u}|^q~\!dz\Big)^{\frac{p}{q}}}+o(1)\Bigg]
\end{equation}
as $h\to\infty$.
It readily follows that $S_{a,b,\gamma}(q)=0$ if $\gamma<b$.

\medskip

The proof of sufficiency is carried out in two steps.

\paragraph{Step 1: purely cylindrical case.}
Here we take  $\gamma=b=0$, hence we are assuming $q\rmH_a>d-k$. In \cite[Corollary 2.1.7/2]{Ma}, it has been proved that $M_{a}(q)>0$ under the additional assumption
$p<d$ (which requires $q\le p^*$). To 
include the case $p\ge d$ and $p<q<\infty$, we argue as in \cite[Appendix B]{MS}, where $p=2=d, a=1$ is assumed.
Put 
$$
\tau=\frac{qp'}{q+p'}~,\qquad \alpha=\frac{q}{\tau}\rmH_a-d+1.
$$
Direct computations show that 
$$
k+\alpha>0~,\qquad 1<\tau\le \frac{d}{d-1}~,\qquad \tau(d-1+\alpha)=q\rmH_a>d-k.
$$
By \cite[Corollary 2.1.7/1]{Ma} we have that there exists $c>0$ such that 
\begin{equation}
\label{eq:Mazya1}
\Big(\int\limits_{\R^d}|y|^{-d+\tau(d-1+\alpha)}|v|^{\tau}~\!dz\Big)^\frac{1}{\tau}\le c\int\limits_{\R^d}|y|^{\alpha}|\nabla v|~\!dz
\quad \text{for any $v\in \mathcal C^1_c(\R^d)$.}
\end{equation}
Fix $u\in \mathcal C^1_c(\R^d)$. 
We use (\ref{eq:Mazya1}) with $v=|u|^{\frac{q}{\tau}}$ (notice that $\frac{q}{\tau}=\frac{q}{p'}+1>1$),
to get
$$
\begin{aligned}
\Big(\int\limits_{\R^d}|y|^{-d+q\rmH_a}|u|^q~\!dz\Big)^\frac{1}{\tau}&=
\Big(\int\limits_{\R^d}|y|^{-d+\tau(d-1+\alpha)}\big||u|^{\frac{q}{\tau}}\big|^{\tau}~\!dz\Big)^\frac{1}{\tau}\\&
\le c\int\limits_{\R^d}|y|^{\alpha}|\nabla |u|^{\frac{q}{\tau}}|~\!dz= c\int\limits_{\R^d}|y|^{\alpha}|u|^\frac{q}{p'}|\nabla u|~\!dz
\\
&\le c \Big(\int\limits_{\R^d}|y|^{a}|\nabla u|^p~\!dz\Big)^\frac1p
\Big(\int\limits_{\R^d}|y|^{(\alpha-\frac{a}{p})p'}|u|^q~\!dz\Big)^\frac{1}{p'}.
\end{aligned}
$$
Since $(\alpha-\frac{a}{p})p'= -d+q\rmH_a$ and $\frac1\tau-\frac{1}{p'}=\frac1q$, 
this proves that $M_{a}(q)=S_{a,0,0}(q)>0$.

\paragraph{Step 2: proof of Theorem \ref{T:Main1} concluded}
To cover  the general case we introduce the parameter 
$$
{t}=a-b+\gamma.
$$
We have ${t}\ge a$ because of assumption $a2)$ in (\ref{eq:assu}), hence $k+{t}>0$. 
Moreover
$q\rmH_t=q\rmH_{a-b+\gamma}>d-k$
by (\ref{eq:CMN_assu2}). Therefore, we can use 
 Step 1 with $a$ replaced by ${t}$  to get  
$$
M_t(q)=\inf_{v\in \mathcal D^{1,p}(\R^d;|y|^{{t}} dz)}
\frac{\displaystyle\int_{\R^d}|y|^{{t}}|\nabla v|^p~\!dz}
{\Big(\displaystyle\int_{\R^d}|y|^{-d+q\rmH_{t}}|v|^q~\!dz\Big)^\frac{p}{q}}>0,
$$
which implies, for any given $u\in \mathcal C^\infty_c(\R^d)$,
\begin{multline*}
\Big(\int\limits_{\R^d}|y|^{-d+q\rmH_{a-b+\gamma}}|z|^{-\gamma\frac{q}{p}} |u|^q~\!dz\Big)^\frac{p}{q}=
\Big(\int\limits_{\R^d}|y|^{-d+q\rmH_t} ||z|^{-\frac{\gamma }{p}}u|^q~\!dz\Big)^\frac{p}{q}\le
c\int\limits_{\R^d}|y|^{{t}}|\nabla (|z|^{-\frac{\gamma}{p}}u)|^p~\!dz\\ 
\le c\int\limits_{\R^d}|y|^{{t}}|z|^{-\gamma}|\nabla u|^p~\!dz+
c \int\limits_{\R^d}|y|^{{t}}|z|^{-\gamma-p}| u|^p~\!dz
\le c
\int\limits_{\R^d}|y|^{{t}}|z|^{-\gamma}|\nabla u|^p~\!dz~\!.
\end{multline*}
Here we used the Hardy--type inequality (\ref{eq:CMN}) with $a,b$ replaced by $t,\gamma$. On the other hand,
since $\gamma\ge b$ by (\ref{eq:assu}), and since $|y||z|^{-1}\le 1$,
we have
$$
\int\limits_{\R^d}|y|^{{t}}|z|^{-\gamma}|\nabla u|^p~\!dz
=\int\limits_{\R^d}|y|^{a}|z|^{-b}\big(|y||z|^{-1}\big)^{\gamma-b }|\nabla u|^p~\!dz
\le \int\limits_{\R^d}|y|^{a}|z|^{-b}|\nabla u|^p~\!dz~\!.
$$
This concludes the proof of the sufficiency.

\medskip

Finally, assume $p<d$. Take a point $z_0=(0,y_0)$, with $|y_0|=1$ and a function
 $u\in \mathcal C^\infty_c(\R^d)\setminus\{0\}$. Then test $S_{a,b,\gamma}(\pstar)$ with $u_h(z)= u((z-z_0)h)$ and  take the limit as $h\to\infty$
to get
$$
S_{a,b,\gamma}(\pstar)\le \frac{\displaystyle\int_{\R^d}|y|^a|z|^{-b}|\nabla u_h|^p}{\Big(\displaystyle\int_{\R^d}|y|^{(a-b+\gamma)\frac{\pstar}{p}}|z|^{-\gamma\frac{\pstar}{p}}|u_h|^\pstar\Big)^\frac{p}{\pstar}}
=
\frac{\displaystyle\int_{\R^d}|\nabla u|^p}{\Big(\displaystyle\int_{\R^d}|u|^\pstar\Big)^\frac{p}{\pstar}}+o(1)~\!,
$$
argue as for (\ref{eq:prima}). This implies  
$$
S_{a,b,\gamma}(\pstar)\le \inf_{u\in \mathcal C^\infty_c(\R^d)} \frac{\displaystyle\int_{\R^d}|\nabla u|^p}{\Big(\displaystyle\int_{\R^d}|u|^\pstar\Big)^\frac{p}{\pstar}}=S~\!.
$$
The theorem is completely proved.
\QED

\section{Proofs of the existence results}
\label{S:Main23}

We start with few preliminaries.

\begin{Lemma}
\label{L:comp}
Let (\ref{eq:CMN_assu}) hold. 
Then 
$\mathcal D^{1,p}(\R^d;|y|^a|z|^{-b}dz)$ is compactly embedded in $L^p_{\rm loc}(\R^d;|y|^a|z|^{-b}dz)$.
\end{Lemma}

\proof
Take a bounded domain $\Omega\subset \R^d$. For  $\eps>0$ we put $\Omega_\eps=\Omega\cap \{(x,y)~|~|y|<\eps\}$.
Then we take $R>0$ so that $\Omega\subset B_R$ and we choose a parameter $t$ such that $0<t<\min\{p,k+a\}$.
For $u\in \mathcal D^{1,p}(\R^d;|y|^a|z|^{-b}dz)$ we estimate
$$
\int\limits_{\Omega_\eps}|y|^a|z|^{-b}|u|^p~\!dz=\int\limits_{\Omega_\eps}(|y|^{a-t}|z|^{-b-p+t} |u|^p)|y|^t|z|^{p-t}~\!dz
\le R^{p-t}\eps^t \ird |y|^{a-t}|z|^{-b-p+t}|u|^p~\!dz.
$$
The Hardy inequality in \cite[Theorem 1.1]{MN_H} applies, and we infer that
\begin{equation}
\label{eq:ball}
\int\limits_{\Omega_\eps}|y|^a|z|^{-b}|u|^p~\!dz\le c \eps^t \ird |y|^{a}|z|^{-b}|\nabla u|^p~\!dz,
\end{equation}
where $c>0$ depends on $\Omega$ but not on $u$.

By the Rellich theorem, we have that for any $\eps>0$ the restriction operator 
$$\mathcal D^{1,p}(\R^d;|y|^a|z|^{-b}dz)\longrightarrow 
L^p(\Omega\setminus\Omega_\eps;|y|^a|z|^{-b}dz)
$$
is compact. Therefore, by 
(\ref{eq:ball}) we have  that the   restriction
$\mathcal D^{1,p}(\R^d;|y|^a|z|^{-b}dz)\longrightarrow 
L^p(\Omega;|y|^a|z|^{-b}dz)$
can be approximated in norm by compact operators. Thus it is compact itself.
This completes the proof of the lemma.
\QED

We denote by 
$(\mathcal D^{1,p}(\R^d;|y|^a|z|^{-b}dz))'$ the dual space of $\mathcal D^{1,p}(\R^d;|y|^a|z|^{-b}dz)$.
Let $\mathcal L_{a,b}$ be the differential operator in (\ref{eq:operator}). 
It holds that 
$\mathcal L_{a,b}(u)\in (\mathcal D^{1,p}(\R^d;|y|^a|z|^{-b}dz))'$,
and
$$
\langle \mathcal L_{a,b}(u), u\rangle= \ird |y|^a|z|^{-b}|\nabla u|^{p}~\!dz
\qquad \text{for  any  $u\in \mathcal D^{1,p}(\R^d;|y|^a|z|^{-b}dz)$.}
$$

\begin{Lemma}
\label{L:operator1}
Assume that (\ref{eq:CMN_assu}) holds. Let $u_h\to 0$ weakly in $\mathcal D^{1,p}(\R^d;|y|^a|z|^{-b}dz)$, and let $\f\in \mathcal C^\infty_c(\R^d)$ be
a (nonnegative) cut-off function. Then 
\begin{itemize}
\item[$i)$] 
$\langle \mathcal L_{a,b}(u_h), \f^pu_h\rangle= \displaystyle\ird |y|^a|z|^{-b}|\nabla (\f u_h)|^{p}~\!dz+o(1)$;
\item[$ii)$] If $p<d$ and $\text{supp}(\f)\subset \R^d\setminus\Sigma_0$,\footnote{Recall that $\Sigma_0$ is the singular set defined in (\ref{eq:Sigma0}).} then
$$
\langle \mathcal L_{a,b}(u_h), \f^pu_h\rangle\ge S\Big(\ird |y|^{a\frac{\pstar}{p}}|z|^{-b\frac{\pstar}{p}}|\f u_h|^{\pstar}~\!dz\Big)^\frac{p}{\pstar}+o(1).
$$
\end{itemize}
\end{Lemma}

\proof
By Lemma \ref{L:comp}
we have that 
$|u_h\nabla \f|\to 0$ in $L^p(\R^d;|y|^a|z|^{-b})$. Therefore,  
$$
\int\limits_{\R^d}|y|^a|z|^{-b}|u_h\nabla \f|^p~\!dz=o(1),
\quad
\int\limits_{\R^d}|y|^a|z|^{-b}|\nabla u_h|^{p-2}(\nabla u_h\cdot\nabla\f^p)u_h~\!dz=o(1)
$$
(the second relation follows from the first one, via H\"older inequality).
Thus we have 
$$
\begin{aligned}
\langle \mathcal L_{a,b}(u_h), \f^pu_h\rangle
&=\int\limits_{\R^d}|y|^a|z|^{-b}|\nabla u_h|^{p-2}\nabla u_h\cdot\nabla(\f^pu_h)~\!dz=
\int\limits_{\R^d}|y|^a|z|^{-b}|\nabla u_h|^{p}\f^{p}~\!dz+o(1)\\
&=\int\limits_{\R^d}|y|^a|z|^{-b}|\nabla (\f u_h)-u_h\nabla \f|^{p}~\!dz+o(1)=\int\limits_{\R^d}|y|^a|z|^{-b}|\nabla (\f u_h)|^{p}~\!dz+o(1),
\end{aligned}
$$
and $i)$ is proved.
\medskip

Now assume  $p<d$ and let $\f\in \mathcal C^\infty_c(\R^d\setminus\Sigma_0)$. Then $|y|^\frac{a}{p}|z|^{-\frac{b}{p}}$ is a smooth function in a neighbourhood of the support of $\f$. Since 
$u_h\to 0$ in $L^p_{\rm loc}(\R^d\setminus\Sigma_0)$ by Rellich theorem, we have that 
$\big|\f u_h\nabla\big(|y|^\frac{a}{p}|z|^{-\frac{b}{p}})\big|\to 0$ in $L^p(\R^d)$, hence
$$
\begin{aligned}
\ird |y|^a|z|^{-b}|\nabla (\f u_h)|^{p}~\!dz &=\int\limits_{\R^d}\big|\nabla \big(|y|^\frac{a}{p}|z|^{-\frac{b}{p}}\f u_h\big)
- \f u_h\nabla\big(|y|^\frac{a}{p}|z|^{-\frac{b}{p}})\big|^{p}~\!dz \\&
=\int\limits_{\R^d}|\nabla \big(|y|^\frac{a}{p}|z|^{-\frac{b}{p}}\f u_h\big)|^{p}~\!dz +o(1).
\end{aligned}
$$
Therefore, the Sobolev inequality gives
$$
\ird |y|^a|z|^{-b}|\nabla (\f u_h)|^{p}~\!dz \ge 
S\Big(\ird |y|^{a\frac{\pstar}{p}}|z|^{-b\frac{\pstar}{p}}|\f u_h|^{\pstar}~\!dz\Big)^\frac{p}{\pstar}+o(1),
$$
which ends the proof, thanks to $i)$.
\QED

\begin{Lemma}
\label{L:localization}
Assume that (\ref{eq:CMN_assu}), (\ref{eq:CMN_assu2}) and (\ref{eq:assu}) hold. Let the sequence $u_h\in 
\mathcal D^{1,p}(\R^d;|y|^a|z|^{-b}dz)$ satisfy 
\begin{equation}
\label{eq:EL}
\mathcal L_{a,b}(u_h)
=|y|^{-d+q\rmH_{a-b+\gamma}}|z|^{-\gamma\frac{q}{p}}|u_h|^{q-2}u_h+o(1)\qquad \text{in $(\mathcal D^{1,p}(\R^d;|y|^a|z|^{-b}dz))'.$}
\end{equation}
If $u_h\to 0$ weakly in $\mathcal D^{1,p}(\R^d;|y|^a|z|^{-b}dz)$, then for any given cut-off function  $\f\in 
\mathcal C^\infty_c(\R^d)$
we have 
\begin{equation}
\label{eq:localization}
\limsup_{h\to\infty}\Big(S_{a,b,\gamma}(q)^\frac{q}{q-p}-\int\limits_{\text{supp}\f}|y|^{-d+q\rmH_{a-b+\gamma}}|z|^{-\gamma\frac{q}{p}}| u_h|^{q}~\!dz\Big)
\Big(\int\limits_{\R^d}|y|^{a}|z|^{-b}|\nabla(\f u_h)|^p~\!dz\Big)^{\frac{q}{q-p}}\le 0.
\end{equation}
\end{Lemma}

\proof
We test (\ref{eq:EL}) with the function $\f^pu_h$. Taking $i)$ of Lemma \ref{L:operator1} into account we get
$$
\begin{aligned}
\int\limits_{\R^d}|y|^a|z|^{-b}&|\nabla (\f u_h)|^{p}~\!dz= \langle \mathcal L_{a,b}(u_h), \f^pu_h\rangle+o(1)
= \int\limits_{\R^d}|y|^{-d+q\rmH_{a-b+\gamma}}|z|^{-\gamma\frac{q}{p}}|u_h|^{q-p}|\f u_h|^p~\!dz+o(1)\\
&\le
\Big(\int\limits_{\text{supp}\f}|y|^{-d+q\rmH_{a-b+\gamma}}|z|^{-\gamma\frac{q}{p}}|u_h|^{q}~\!dz\Big)^{\frac{q-p}{q}}
\Big(\int\limits_{\R^d}|y|^{-d+q\rmH_{a-b+\gamma}}|z|^{-\gamma\frac{q}{p}}|\f u_h|^q~\!dz\Big)^{\frac{p}{q}}+o(1)\\
&
\le \frac{1}{S_{a,b,\gamma}(q)} \Big(\int\limits_{\text{supp}\f}|y|^{-d+q\rmH_{a-b+\gamma}}|z|^{-\gamma\frac{q}{p}}|u_h|^{q}~\!dz\Big)^{\frac{q-p}{q}}
\int\limits_{\R^d}|y|^a|z|^{-b}|\nabla (\f u_h)|^{p}~\!dz+o(1)~\!,
\end{aligned}
$$
which easily implies (\ref{eq:localization}).
\QED

\begin{Lemma}
\label{L:operator2}
Assume that (\ref{eq:assu_cyl}) is satisfied. Let $b<p\rmH_a$, and  $v\in \mathcal D^{1,p}(\R^d;|y|^a|z|^{-b}dz)$. Then
$$
\Big(\int\limits_{\R^d}|y|^a|z|^{-b}|\nabla v|^p~\!dz\Big)^\frac1p+\tfrac{|b|}{p}\Big(\int\limits_{\R^d}|y|^a|z|^{-b-p}|v|^p~\!dz\Big)^\frac1p
\ge M_{a}(q)^\frac1p \Big(\int\limits_{\R^d}|y|^{-d+q\rmH_a}|z|^{-b\frac{q}{p}}|v|^q~\!dz\Big)^\frac{1}{q}.
$$
\end{Lemma}

\proof
For $v\in \mathcal C^\infty_c(\R^d\setminus\{0\})$ use the triangle inequality and then the Maz'ya inequality in (\ref{eq:Ma})
to get
$$
\begin{aligned}
\Big(\int\limits_{\R^d}|y|^a|z|^{-b}|\nabla v|^p~\!dz\Big)^\frac1p&= 
\Big(\int\limits_{\R^d}|y|^a\big|\nabla (|z|^{-\frac{b}{p}}v)+\tfrac{b}{p}|z|^{-\frac{b}{p}-2}z v\big|^p~\!dz\Big)^\frac1p
\\& \ge 
\Big(\int\limits_{\R^d}|y|^a|\nabla (|z|^{-\frac{b}{p}}v)|^p~\!dz\Big)^\frac1p- 
\tfrac{|b|}{p}\Big(\int\limits_{\R^d}|y|^a|z|^{-b-p}|v|^p~\!dz\Big)^\frac1p\\
&\ge M_{a}(q)^\frac1p\Big(\int\limits_{\R^d}|y|^{-d+q\rmH_a}||z|^{-\frac{b}{p}}v|^q~\!dz\Big)^\frac{1}{q}
-\tfrac{|b|}{p}\Big(\int\limits_{\R^d}|y|^a|z|^{-b-p}|v|^p~\!dz\Big)^\frac1p.
\end{aligned}
$$
The conclusion  follows, since $\mathcal C^\infty_c(\R^d\setminus\{0\})$ is dense in $\mathcal D^{1,p}(\R^d;|y|^a|z|^{-b}dz)$
by \cite[Lemma 2.1]{CMN}. 
\QED

We conclude this preliminary section with a lemma of independent interest concerning the weak limit of sequences 
satisfying (\ref{eq:EL}). The argument is essentially due to Evans \cite{Ev}, see also \cite{Sa}.

\begin{Lemma}
\label{L:PSSa}
Assume that (\ref{eq:CMN_assu}), (\ref{eq:CMN_assu2}) and (\ref{eq:assu}) hold. Let $u_h\in 
\mathcal D^{1,p}(\R^d;|y|^a|z|^{-b}dz)$ be a sequence   satisfying  (\ref{eq:EL}).
If $u_h\to u$ weakly in $\mathcal D^{1,p}(\R^d;|y|^a|z|^{-b}dz)$, then $u$ solves 
\begin{equation}
\label{eq:ELu}
\mathcal L_{a,b}(u)
=|y|^{-d+q\rmH_{a-b+\gamma}}|z|^{-\gamma\frac{q}{p}}|u|^{q-2}u.
\end{equation}
\end{Lemma}

\proof
The case $p=2$ is straightforward. To deal with the nonlinear case $p\neq 2$ 
we adapt the argument in \cite[Step 1.2]{Sa}. 

In this proof we omit the integration element $dz$ and put
$$
\nu(z)= |y|^{-d+q\rmH_{a-b+\gamma}}|z|^{-\gamma\frac{q}{p}}.
$$
Since $u_h\to u$ in $L^p_{\rm loc}(\R^d\setminus\Sigma_0)$ by Rellich Theorem, we can assume that $u_h\to u$ a.e. on $\R^d$.
Further, $u_h$ is bounded in $L^q(\R^d;\nu(z)dz)$ by Theorem \ref{T:Main1}, thus 
$|u_h|^{q-2}u_h\to |u|^{q-2}u$ 
weakly in $L^{q'}(\R^d;\nu(z)dz)$. 
In particular, for any fixed test function $\eta\in \mathcal C^\infty_c(\R^d)$ we have that 
$$
\ird\nu(z)|u_h|^{q-2}u_h\eta= \ird\nu(z)|u|^{q-2}u\eta+o(1).
$$
We  need to show that $\langle \mathcal L_{a,b}(u_h),\eta\rangle=\langle \mathcal L_{a,b}(u),\eta\rangle+o(1)$, that is,
\begin{equation}
\label{eq:need}
\ird|y|^a|z|^{-b}F_h
\cdot\nabla \eta= \ird|y|^a|z|^{-b}F
\cdot\nabla \eta+o(1)~\!,
\end{equation}
where
$$
F_h= |\nabla u_h|^{p-2}\nabla u_h~,\qquad F= |\nabla u|^{p-2}\nabla u~\!.
$$
To conclude the proof we will show that, up to a subsequence,
\begin{equation}
\label{eq:need1}
F_h\to F\qquad \text{a.e. in $\R^d$.}
\end{equation}
In fact, since the sequence of vectorfields $F_h$ is bounded in $L^{p'}(\R^d;|y|^a|z|^{-b}dz)^d$, then   (\ref{eq:need1}) implies that $F_h\to F$ weakly in $L^{p'}(\R^d;|y|^a|z|^{-b}dz)^d$, and 
(\ref{eq:need}) follows.

\medskip
To prove (\ref{eq:need1}) fix  $R, \delta >0$. By  the Egorov theorem, there exists $E_\delta \subset B_R=\{|z|<R\}$ such that 
\[
\int\limits_{B_R \setminus E_\delta}|y|^a |z|^{-b} dz < \delta\quad\text{and}\quad \text{$u_h \to u$ uniformly in $E_\delta$.}
\] 
Next, fix $\eps>0$. For $h$  large, we have that
\begin{equation}\label{eq:unifconv}
|u_h - u | <\eps \quad \text{ in }E_\delta\,. 
\end{equation}
Following \cite{Sa}, we truncate the sequence $u_h-u$ as follows:
$$
\beta^\eps_h=\begin{cases}
u_h-u&\text{if $|u_h-u|\le \eps$}\\
\pm \eps&\text{if $\pm(u_h-u)>\eps$.}
\end{cases}
$$
Notice that 
$$|\beta^\eps_h|\le \min\{\eps,|u_h-u|\}~\!,
$$
thus $\beta^\eps_h\to 0$ a.e. on $\R^d$ as $h\to \infty$. In addition, by
using the inequality $|\nabla \beta^\eps_h|\le |\nabla(u_h-u)|$, one can plainly prove that $(\beta^\eps_h)_h$ is 
a bounded sequence in $\mathcal D^{1,p}(\R^d;|y|^a|z|^{-b}dz)$, and that
\begin{equation}
\label{eq:fibeta}
\beta^\eps_h\to 0\qquad \text{weakly in $\mathcal D^{1,p}(\R^d;|y|^a|z|^{-b}dz)$.}
\end{equation}

Next, take a cut-off function $\f\in \mathcal C^\infty_c(\R^d)$ such that $0\le \f\le 1$ and $\f\equiv 1$ on $B_R$. 
Then 
\begin{equation}
\label{eq:fibeta2}
\ird |y|^a|z|^{-b} F\cdot \nabla (\f\beta^\eps_h)=o(1)\qquad \text{as $h\to\infty$}~\!,
\end{equation}
because $\varphi \beta^\eps_h\to 0$ weakly in $\mathcal D^{1,p}(\R^d;|y|^a|z|^{-b}dz)$. Thanks to (\ref{eq:fibeta2}), (\ref{eq:fibeta}) and Lemma \ref{L:comp}, we infer that
$$
\langle \mathcal L_{a,b}(u_h),\varphi \beta^\eps_h\rangle 
=\ird |y|^a|z|^{-b} (F_h-F)\cdot \nabla(\f\beta^\eps_h)+o(1)
=\ird |y|^a|z|^{-b} \big((F_h-F)\cdot \nabla \beta^\eps_h\big)~\!\f+o(1)~\!.
$$
Now we use the elementary but crucial convexity inequality
\begin{equation}\label{eq:convexity}
(F_h - F)\cdot (\nabla u_h - \nabla u ) > 0\qquad \text{a.e. on $\{\nabla u_h\neq \nabla u\}$}\,.
\end{equation}
This  gives $(F_h - F)\cdot \nabla \beta^\eps_h\ge 0$ a.e. on $\R^d$. In addition, 
thanks to (\ref{eq:unifconv}), we have that  $\nabla\beta^\eps_h=\nabla u_h-\nabla u$ on $E_\delta$.
Thus (\ref{eq:convexity}) gives 
$$
\langle \mathcal L_{a,b}(u_h),\varphi \beta^\eps_h\rangle\ge \int\limits_{E_\delta} |y|^a|z|^{-b} (F_h-F)\cdot (\nabla u_h-\nabla u)+o(1),
$$
because  $\f\equiv 1$ on $E_\delta$.
Therefore, by testing (\ref{eq:EL}) with $\f\beta^\eps_h$ we get
$$
0\le
\int\limits_{E_\delta} |y|^a|z|^{-b} (F_h-F)\cdot (\nabla u_h-\nabla u)\le 
\ird\nu(z)|u_h|^{q-2}u_h(\f\beta^\eps_h)+o(1)\le
\eps\int\limits_{\text{supp}\f}\nu(z)|u_h|^{q-1}+o(1)~\!.
$$
Since in addition $\nu(z)\in L^1_{\rm loc}(\R^d)$, we can apply the H\"older inequality to obtain the estimate 
\[
0 \leq \limsup_{h \to \infty} \int_{E_\delta} |y|^a |z|^{-b}(F_h - F) \cdot (\nabla u_h - \nabla u ) \leq c \eps,
\]
with $c>0$ not depending on  the arbitrarily chosen parameter $\eps>0$. Taking (\ref{eq:convexity}) and Fatou's Lemma into account, we infer that
$\nabla u_h - \nabla u\to 0$ a.e. in $E_\delta$. Then (\ref{eq:need1}) follows, since $\delta, R>0$ were arbitrarily chosen.
This concludes the proof of  the Lemma.
\QED

\subsection{Proof of Theorem \ref{T:Main2}}
We start by constructing a well-behaved minimizing sequence for the minimization problem (\ref{eq:HS_constant}). 

Let us introduce some notation.
For $R>0$, $\ell\in \mathbb N$, $\xi=(\xi_1,\cdots,\xi_\ell)\in\R^\ell$ and $t>0$ we put
\begin{equation}
\label{eq:notation}
B^\ell_t=\big\{\xi\in\R^\ell~|~ |\xi|^2=\sum_{j=1}^\ell \xi_j^2<t^2~\big\}~,\qquad C^\ell_t=\big\{\xi\in\R^\ell~|~ |\xi|_\infty=\max_{j=1,\dots,\ell} |\xi_j|<t~\big\},
\end{equation}
so that $B^\ell_t+\xi_0$ and $C^\ell_t+\xi_0$ are centered at the point $\xi_0\in\R^\ell$.

Furthermore, to simplify the notation we omit the integration element $dz$ and  
put 
$$
S_{a,b,\gamma}=S_{a,b,\gamma}(q)~,\quad \Theta=-d+q\rmH_{a-b+\gamma}=-d+(d-p+a-b+\gamma)\frac{q}{p},
$$
so that
$$
S_{a,b,\gamma}= \inf_{u\in \mathcal D^{1,p}(\R^d;|y|^{a}|z|^{-b} dz)}
\frac{\displaystyle\int_{\R^d}|y|^{a}|z|^{-b}|\nabla u|^p}
{\Big(\displaystyle\int_{\R^d}|y|^{\Theta}|z|^{-\gamma\frac{q}{p}}|u|^q\Big)^\frac{p}{q}}~\!.
$$

We take a minimizing sequence $v_h$, satisfying the following normalization condition
$$
S_{a,b,\gamma}^\frac{q}{q-p}=\ird |y|^{\Theta}|z|^{-\gamma\frac{q}{p}}|v_h|^q=\ird |y|^{a}|z|^{-b}|\nabla v_h|^p+o(1).
$$
By  Ekeland's variational principle \cite{Ek}, we can assume that 
$$
\mathcal L_{a,b}(v_h)
=|y|^{\Theta}|z|^{-\gamma\frac{q}{p}}|v_h|^{q-2}v_h+o(1)\qquad \text{in}~~(\mathcal D^{1,p}(\R^d;|y|^a|z|^{-b}dz))',
$$
where $\mathcal L_{a,b}$ is defined in (\ref{eq:operator}).
By using in a standard way the concentration function
$$
Q_h(t)=\sup_{\xi\in\Sigma_0}\int\limits_{C^d_t+t\xi}|y|^\Theta|z|^{-\gamma\frac{q}{p}}|v_h|^q~,\quad t>0,
$$
one can find sequences $t_h>0$, $\xi_h=(x_h,0)\in\Sigma_0$ such that 
$$
Q_h(t_h)=\frac12 S_{a,b,\gamma}^\frac{q}{q-p}~,\quad 
\int\limits_{C^d_{t_h}+t_h\xi_h}|y|^\Theta|z|^{-\gamma\frac{q}{p}}|v_h|^q\ge \frac14S_{a,b,\gamma}^\frac{q}{q-p}.
$$
It follows that the rescaled sequence
$u_h(z)=t_h^{\rmH_{a-b}}v_h(t_hz)$ is minimizing for $S_{a,b,\gamma}$, satisfies (\ref{eq:EL}) and moreover
\begin{eqnarray}
\label{eq:tu_0}
\ird |y|^{\Theta}|z|^{-\gamma\frac{q}{p}}|u_h|^q&=&
S_{a,b,\gamma}^\frac{q}{q-p}+o(1)~\!,\\
\label{eq:tu_1}
\int\limits_{C^d_1+\xi_h}|y|^\Theta|z|^{-\gamma\frac{q}{p}}|u_h|^q&\ge&\frac14S_{a,b,\gamma}^\frac{q}{q-p}~\!, \\
\label{eq:tu_2}
\int\limits_{C^d_1+\xi}|y|^\Theta|z|^{-\gamma\frac{q}{p}}|u_h|^q&\le&  \frac12 S_{a,b,\gamma}^\frac{q}{q-p}\quad\text{for any $\xi\in\Sigma_0$.}
\end{eqnarray}
Up to a subsequence, we can assume that $u_h\to u$ weakly in $\mathcal D^{1,p}(\R^d;|y|^a|z|^{-b}dz)$. Then $u$ solves 
(\ref{eq:ELu}) by Lemma \ref{L:PSSa} and in particular
\begin{equation}
\label{eq:u}
\ird|y|^a|z|^{-b}|\nabla u|^p=\ird |y|^{\Theta}|z|^{-\gamma\frac{q}{p}}|u|^q.
\end{equation}
To conclude the proof, it remains to show that $u \neq 0$. Indeed, if this is the case,  we can test $S_{a,b,\gamma}$ with $u$
to obtain, via (\ref{eq:u}),
$$
S_{a,b,\gamma} \le 
\left(\int_{\R^d} |y|^{a} |z|^{-b} |\nabla u|^p \right)^{\frac{q-p}{q}} 
\le \liminf_{h \to \infty} \left(\int_{\R^d} |y|^{a} |z|^{-b} |\nabla u_h|^p \right)^{\frac{q-p}{q}} = S_{a,b,\gamma}~\!.
$$
This readily ensures that $u$ achieves $S_{a,b,\gamma}$.

To prove that $u\neq 0$ we first show that the sequence $\xi_h$ is bounded.
Here, the assumption $\gamma>b$ plays a crucal role. We start by taking a new parameter
$\tgamma$ such that
$$
b<\tgamma<\gamma~,\qquad q\rmH_{a-b+\tgamma}>d-k~\!. 
$$
Notice that $S_{a,b,\tgamma}>0$ by Theorem \ref{T:Main1}. 

Assume $|\xi_h|_\infty\ge2$ for some  $h$. Since $\xi_h=(x_h,0)\in \Sigma_0$, then for 
$z=(x,y)\in C^d_1+\xi_h$ it holds that
\begin{equation}
\label{eq:ine_cuve}
|z|_\infty\ge  |\xi_h|_\infty-1\ge \frac{1}{2}|\xi_h|_\infty~,~\qquad |y|_\infty<1.
\end{equation}
It follows that on the cube $C^d_1+\xi_h$ we have
$|y||z|^{-1} \le c |y|_\infty|z|_\infty^{-1} 
\le 
c|\xi_h|^{-1}$, and thus 
$$
|y|^{\Theta}|z|^{-\gamma\frac{q}{p}}
=(|y||z|^{-1})^{(\gamma-\tgamma)\frac{q}{p}}|y|^{-d+q\rmH_{a-b+\tgamma}} |z|^{-\tgamma\frac{q}{p}}
\le
c{|\xi_h|^{-(\gamma-\tgamma)\frac{q}{p}}}|y|^{-d+q\rmH_{a-b+\tgamma}} |z|^{-\tgamma\frac{q}{p}}.
$$
We substitute this estimate in (\ref{eq:tu_1}) to get
$$
\frac14S_{a,b,\gamma}^\frac{q}{q-p}\le
c{|\xi_h|^{-(\gamma-\tgamma)\frac{q}{p}}} \int\limits_{C_1^d+\xi_h}|y|^{-d+q\rmH_{a-b+\tgamma}} |z|^{-\tgamma\frac{q}{p}}|u_h|^q
\le  c{|\xi_h|^{-(\gamma-\tgamma)\frac{q}{p}}} \Big(\frac{1}{S_{a,b,\tgamma}}
\int\limits_{\R^d}|y|^{a} |z|^{-b}|\nabla u_h|^p\Big)^\frac{q}{p}.
$$
This gives the boundedness of the sequence $\xi_h\in\Sigma_0$, because $u_h$ is bounded in $\mathcal D^{1,p}(\R^d;|y|^a|z|^{-b}dz)$.
Thus we can take $R>0$ large enough, so that $C^d_1+\xi_h\subset C_R^{d-k}\times C^k_1$
for any $h\ge 1$. The inequality (\ref{eq:tu_1}) then gives
\begin{equation}
\label{eq:included}
\int\limits_{C_R^{d-k}\times C^k_1} |y|^\Theta|z|^{-\gamma\frac{q}{p}}|u_h|^q\ge \frac14S_{a,b,\gamma}^\frac{q}{q-p}.
\end{equation}

Assume by contradiction  that $u_h$ converges weakly to $0$. 
Take any cut-off function $\f$ having support in a cube of side $1$ about a point $\xi\in\Sigma_0$.
Then (\ref{eq:tu_2}) and (\ref{eq:localization}) give
$$
\int\limits_{\R^d}|y|^{a}|z|^{-b}|\nabla(\f u_h)|^p~\!dz=o(1),
$$
and therefore $\f u_h\to 0$ in $L^q(\R^d;|y|^\Theta|z|^{-\gamma\frac{q}{p}} dz)$. This implies that 
$u_h\to 0$ in $L^q({C^d_\frac12}+\xi;|y|^\Theta|z|^{-\gamma\frac{q}{p}} dz)$ for any  $\xi \in \Sigma_0$.
Finally, a compactness argument gives
$$
\int\limits_{C^{d-k}_R\times C^k_\frac12}|y|^\Theta|z|^{-\gamma\frac{q}{p}}|u_h|^q=o(1),
$$
which compared with (\ref{eq:included}) implies 
\begin{equation}
\label{eq:KK}
\int\limits_{A_R}|y|^\Theta|z|^{-\gamma\frac{q}{p}}|u_h|^q\ge \frac14S_{a,b,\gamma}^\frac{q}{q-p}+o(1)~,
\end{equation}
where
\begin{equation}
\label{eq:Omega}
A_R:=C^{d-k}_R\times(C^k_1\setminus \overline{C^k_\frac12})=\big\{(x,y)\in\R^d~|~|x|_\infty<R~,~\frac12<|y|_\infty<1~\!\big\}
\Subset \R^d\setminus\Sigma_0.
\end{equation}

\medskip

\paragraph{Proof of $i)$}
If $p\ge d$, or if $p<d$ and $q<p^*$, then 
$u_h\to 0$ strongly in $L^q(A_R)$ by the Rellich theorem. 
This contradicts (\ref{eq:KK}) and concludes the proof in this case.

\paragraph{Proof of $ii)$}
Now let $q=p^*$ and assume $S_{a,b,\gamma}(\pstar)<S$. Recall that $u_h$ satisfies (\ref{eq:EL}) and (\ref{eq:tu_0}) which, in the critical case, transform into 
\begin{gather}
\label{eq:EL*}
\mathcal L_{a,b}(u_h)
=|y|^{\Theta}|z|^{-\gamma\frac{\pstar}{p}}|u_h|^{p^*-2}u_h+o(1)~,\quad\Theta=(a-b+\gamma)\frac{\pstar}{p}~\!,\\
\label{eq:mass}
\int\limits_{\R^d}|y|^\Theta|z|^{-\gamma\frac{\pstar}{p}}|u_h|^{\pstar}= 
S_{a,b,\gamma}^\frac{\pstar}{\pstar-p}+o(1).
\end{gather}
Take a cut-off function $\f\in \mathcal C^\infty_c(\R^d\setminus\Sigma_0)$ such that $\f\equiv 1$ on the set
${A_R}$ in (\ref{eq:Omega}). 
We estimate 
$$
|y|^\Theta|z|^{-\gamma\frac{\pstar}{p}}= |y|^{a\frac{\pstar}{p}}|z|^{-b\frac{\pstar}{p}}(|y||z|^{-1})^{(\gamma-b)\frac{\pstar}{p}}
\le |y|^{a\frac{\pstar}{p}}|z|^{-b\frac{\pstar}{p}}
$$
because $\gamma\ge b$ and $|y|\le |z|$. Combining this information with $ii)$ in Lemma \ref{L:operator1}  we get
$$
S\Big(\ird|y|^\Theta|z|^{-\gamma\frac{\pstar}{p}}|\f u_h|^{\pstar}\Big)^\frac{p}{\pstar}
\le 
S\Big(\int|y|^{a\frac{\pstar}{p}}|z|^{-b\frac{\pstar}{p}}|\f u_h|^{\pstar}\Big)^\frac{p}{\pstar} \le 
\langle \mathcal L_{a,b}(u_h), \f^pu_h\rangle +o(1).
$$
Thus (\ref{eq:EL*}),  H\"older inequality and (\ref{eq:mass}) give
$$
\begin{aligned}
S\Big(\ird|y|^\Theta|z|^{-\gamma\frac{\pstar}{p}}|\f u_h|^{\pstar}\Big)^\frac{p}{\pstar}
&\le
\int\limits_{\R^d}|y|^\Theta|z|^{-\gamma\frac{\pstar}{p}}|u_h|^{\pstar-p}|\f u_h|^p+o(1)\\
&\le \Big(\int\limits_{\R^d}|y|^\Theta|z|^{-\gamma\frac{\pstar}{p}}|u_h|^{\pstar}\Big)^\frac{\pstar-p}{\pstar}
\Big(\ird|y|^\Theta|z|^{-\gamma\frac{\pstar}{p}}|\f u_h|^{\pstar}\Big)^\frac{p}{\pstar}+o(1)\\ &
=S_{a,b,\gamma}~\!\Big(\ird|y|^\Theta|z|^{-\gamma\frac{\pstar}{p}}|\f u_h|^{\pstar}\Big)^\frac{p}{\pstar}+o(1).
\end{aligned}
$$
Since $S_{a,b,\gamma}<S$ and $\f\equiv 1$ on $A_R$, we infer
$$
(S- S_{a,b,\gamma})
\Big(\int\limits_{{A_R}}|y|^\Theta|z|^{-\gamma\frac{\pstar}{p}}|u_h|^{\pstar}\Big)^\frac{p}{\pstar}
\le (S- S_{a,b,\gamma})
\Big(\ird|y|^\Theta|z|^{-\gamma\frac{\pstar}{p}}|\f u_h|^{\pstar}\Big)^\frac{p}{\pstar}
\le o(1).
$$
This contradicts (\ref{eq:included}) and concludes the proof of the theorem.
\QED

\subsection{Proof of Theorem \ref{T:Main3}}

We use again the notation (\ref{eq:notation}) for balls and cubes in $\R^d$, and omit the integration element $dz$. Also, we 
put $M_{a}:=M_{a}(q)$ and $S_{a,b}:=S_{a,b,b}(q)$, compare with (\ref{eq:Ma}) and (\ref{eq:Mb}), so that
$$
M_{a}
= \inf_{u\in \mathcal D^{1,p}(\R^d;|y|^{a} dz)}
\frac{\displaystyle\int_{\R^d}|y|^{a}|\nabla u|^p}
{\Big(\displaystyle\int_{\R^d}|y|^{-d+q\rmH_a}|u|^q\Big)^\frac{p}{q}}~\!,
\quad
S_{a,b}=\inf_{u\in \mathcal D^{1,p}(\R^d;|y|^{a}|z|^{-b} dz)}
\frac{\displaystyle\int_{\R^d}|y|^{a}|z|^{-b}|\nabla u|^p}
{\Big(\displaystyle\int_{\R^d}|y|^{-d+q\rmH_a}|z|^{-b\frac{q}{p}}|u|^q\Big)^\frac{p}{q}}~\!.
$$

\paragraph{Proof of $i)$}
Take $x_0\in\R^{d-k}$ such that $|x_0|=1$ and put $z_0=(x_0,0)$. Take a nontrivial function ${u}\in \mathcal C^\infty_c(\R^d)$ and put ${u}_h(z)={u}(z-hz_0)$. 
By (\ref{eq:prima}) we have that
$$
S_{a,b}\le \frac{\displaystyle \int_{\R^d}|y|^{a}|z|^{-b}|\nabla {u}_h|^p}
{\Big(\displaystyle \int_{\R^d}|y|^{-d+q\rmH_a}|z|^{-b\frac{q}{p}}|{u}_h|^q\Big)^{\frac{p}{q}}}
=\frac{\displaystyle \int_{\R^d}|y|^{a}|\nabla {u}|^p}
{\Big(\displaystyle \int_{\R^d}|y|^{-d+q\rmH_a}|{u}|^q\Big)^{\frac{p}{q}}}+o(1)
$$
as $h\to \infty$. Since ${u}$ was arbitrarily chosen, we infer that
$$
S_{a,b}\le\inf_{{u}\in {\mathcal C}^\infty_c(\R^d)} 
\frac{\displaystyle \int_{\R^d}|y|^{a}|\nabla {u}|^p}
{\Big(\displaystyle \int_{\R^d}|y|^{-d+q\rmH_a}|{u}|^q\Big)^{\frac{p}{q}}}=M_{a},
$$
which proves $i)$.

\paragraph{Proof of $ii)$} 
Recall that now we assume $S_{a,b}<M_a$.
As in the proof of Theorem \ref{T:Main2}, we apply Ekeland's variational principle \cite{Ek} to construct a bounded minimizing sequence $v_h$
which satisfies also 
$$
\begin{gathered}
\langle \mathcal L_{a,b}(v_h), v_h \rangle
=|y|^{-d+q\rmH_a}|z|^{-b\frac{q}{p}}|v_h|^{q-2}v_h+o(1)
\quad \text{in $(\mathcal D^{1,p}(\R^d;|y|^a|z|^{-b}dz))'$.}
\\
\ird |y|^{-d+q\rmH_a}|z|^{-b\frac{q}{p}}|v_h|^q=\ird |y|^{a}|z|^{-b}|\nabla v_h|^p+o(1)=S_{a,b}^\frac{q}{q-p}+o(1).
\end{gathered}
$$
Thanks to Lemma \ref{L:operator2} we have
$$
\begin{aligned}
\liminf_{h\to \infty}\tfrac{|b|}{p}\Big(\int\limits_{\R^d}|y|^{a}|z|^{-b-p}|v_h|^p\Big)^\frac1p\ge 
M_{a}^\frac1p S_{a,b}^\frac{1}{q-p}-S_{a,b}^\frac{q}{p(q-p)}\ge 
S_{a,b}^\frac{1}{q-p}\big(M_{a}^\frac{1}{p}-S_{a,b}^\frac{1}{p}\big)>0.
\end{aligned}
$$
We now need to properly rescale the sequence $v_h$ to obtain a new minimizing sequence
which can not concentrate at points at infinity corresponding to the directions of $\Sigma_0$. To this end, 
the concentration function $Q_h(t)$ we used 
in the proof of Theorem \ref{T:Main2} has to be replaced by
$$
\mathcal  Q_h(t)=\sup_{\xi\in\Sigma_0}\int\limits_{C_t^d+t\xi}|y|^a|z|^{-b-p}|v_h|^p~,\quad t>0~\!.
$$
We can find sequences $t_h>0$ and $\xi_h=(x_h,0)$ in the singular set $\Sigma_0$, such that the rescaled sequence
$u_h(z)=t_h^{\rmH_{a-b}}v_h(t_hz)$ is minimizing for $S_{a,b}$, satisfies (\ref{eq:EL}) and moreover
\begin{eqnarray}
\label{eq:uh_energy}
\ird |y|^{-d+q\rmH_a}|z|^{-b\frac{q}{p}}|u_h|^q&=&
S_{a,b}^\frac{q}{q-p}+o(1),\\
\label{eq:tu_b1}
\int\limits_{C^d_1+\xi_h}|y|^a|z|^{-b-p}|u_h|^p&\ge& \delta\\
\label{eq:tu_b2}
\sup_{\xi\in\Sigma_0} \int\limits_{C^d_1+\xi}|y|^a|z|^{-b-p}|u_h|^p&\le& 2\delta
\end{eqnarray}
for some small $\delta>0$ to be chosen later.

Up to a subsequence, we can assume that $u_h\to u$ weakly in $\mathcal D^{1,p}(\R^d;|y|^a|z|^{-b}dz)$.

The next step consists in showing that $u\neq 0$, from which the conclusion follows by the same argument as in the proof of Theorem~\ref{T:Main2}.

\medskip
Firstly, we need to prove that the sequence $\xi_h$ is bounded.  
Assume that $|\xi_h|_\infty\ge2$ for some  $h$. Since the assumption $q(d-p)\le dp$  implies 
$-d+\frac{pq}{q-p}\ge 0$ and since $|y|_\infty\le 1$ on $C_1^d+\xi_h$, then the relations in (\ref{eq:ine_cuve}) give
$$
|y|^{-d+\frac{pq}{q-p}}|z|^{-\frac{pq}{q-p}}\le c|z|^{-\frac{pq}{q-p}}
\le  c|\xi_h|^{-\frac{pq}{q-p}}\qquad \text{on $C^d_1+\xi_h$.}
$$
Thus from (\ref{eq:tu_b1}) it follows that 
$$
\begin{aligned}
\delta\le \int\limits_{C^d_1+\xi_h} |y|^a|z|^{-b-p}|u_h|^p&\le    
 \Big( \int\limits_{C^d_1+\xi_h} |y|^{-d+q\rmH_a}|z|^{-b\frac{q}{p}}|u_h|^q\Big)^\frac{p}{q}
  \Big(\int\limits_{C^d_1+\xi_h} |y|^{-d+\frac{pq}{q-p}}|z|^{-\frac{pq}{q-p}}\Big)^\frac{q-p}{q}
  \\&
  \le c\Big( \ird |y|^{-d+q\rmH_a}|z|^{-b\frac{q}{p}}|u_h|^q\Big)^\frac{p}{q} |\xi_h|^{-p}
  \end{aligned}
  $$
which, together with (\ref{eq:uh_energy}), gives the boundedness of the sequence $\xi_h$. Thus  we can find $R>0$ such that 
$C^d_1+\xi_h\subset C_R^{d-k}\times C^k_1$
for any $h\ge 1$, and in particular
\begin{equation}
\label{eq:included0}
\int\limits_{C_R^{d-k}\times C^k_1} |y|^a|z|^{-b-p}|u_h|^p\ge \delta\quad\text{for any $h\ge 1$.}
\end{equation}

\medskip
We now assume by contradiction that $u_h\to 0$ weakly in $\mathcal D^{1,p}(\R^d;|y|^a|z|^{-b}dz)$. 
Take any point $\xi\in\Sigma_0$ and any cut-off function $\f\in \mathcal C^\infty_c(C^d_1+\xi)$
such that $\f\equiv 1 $ on $C^d_{\frac34}+\xi$.
We use Lemma \ref{L:operator2} with $v=\f u_h$ to get  
\begin{equation}
\label{eq:compare}
M_{a}^\frac1p \Big(\int\limits_{\R^d}|y|^{-d+q\rmH_a}|z|^{-b\frac{q}{p}}|\f u_h|^q\Big)^\frac{1}{q}
\le
\Big(\int\limits_{\R^d}|y|^a|z|^{-b}|\nabla (\f u_h)|^p\Big)^\frac1p+
\tfrac{|b|}{p}\Big(\int\limits_{\R^d}|y|^a|z|^{-b-p}|\f u_h|^p\Big)^\frac1p.
\end{equation}
Next, Lemma \ref{L:operator1}, (\ref{eq:EL}), H\"older inequality and (\ref{eq:uh_energy}) give 
$$
\begin{aligned}
\int\limits_{\R^d}|y|^a|z|^{-b}|\nabla (\f u_h)|^p&=\langle \mathcal L_{a,b}(u_h), \f^pu_h\rangle +o(1)
=\int\limits_{\R^d}|y|^{-d+q\rmH_a}|z|^{-b\frac{q}{p}}|u_h|^{q-p}|\f u_h|^p +o(1)\\ &
\le
\Big(\ird |y|^{-d+q\rmH_a}|z|^{-b\frac{q}{p}}|u_h|^{q}~\!\Big)^{\frac{q-p}{q}}
\Big(\int\limits_{\R^d}|y|^{-d+q\rmH_a}|z|^{-b\frac{q}{p}}|\f u_h|^q~\!\Big)^{\frac{p}{q}}+o(1)\\
&= S_{a,b}~\!\Big(\int\limits_{\R^d}|y|^{-d+q\rmH_a}|z|^{-b\frac{q}{p}}|\f u_h|^q~\!\Big)^{\frac{p}{q}}+o(1).
\end{aligned}
$$
Comparing with (\ref{eq:compare}) we infer 
$$
\big(M_{a}^\frac{1}{p}-S_{a,b}^\frac{1}{p}\big)\Big(\int\limits_{\R^d}|y|^{-d+q\rmH_a} |z|^{-\frac{b}{p}}|\f u_h|^{q}\Big)^\frac{1}{q}\le 
\tfrac{|b|}{p}\Big(\int\limits_{\R^d}|y|^a|z|^{-b-p}|\f u_h|^{p}\Big)^\frac{1}{p}+o(1).
$$
Since $M_{a}>S_{a,b}$, and since the cut-off function $\f\in \mathcal C^\infty_c(C^d_1+\xi)$ satisfies $\f\equiv 1$ on $C^d_\frac34+\xi$, we get
$$
\begin{aligned}
\big(M_{a}^\frac{1}{p}-S_{a,b}^\frac{1}{p}\big)\Big(\int\limits_{C^d_\frac34+\xi}&|y|^{-d+q\rmH_a} |z|^{-\frac{b}{p}}|u_h|^{q}\Big)^\frac{1}{q}
\le \tfrac{|b|}{p}\Big(\int\limits_{C^d_1+\xi}|y|^a|z|^{-b-p}|u_h|^{p}\Big)^\frac{1}{p}+o(1)\le \tfrac{|b|}{p}(2\delta)^\frac{1}{p}+o(1).
\end{aligned}
$$
by (\ref{eq:tu_b2}). So, for $\delta$ small enough and $h$ large we obtain that
\begin{equation}
\label{eq:small}
\int\limits_{C^d_{\frac34}+\xi}|y|^{-d+q\rmH_a} |z|^{-b\frac{q}{p}}|u_h|^{q}\le 
\frac12 S_{a,b}^\frac{q}{q-p}.
\end{equation}
Now,  take another cut-off function $\tilde\f\in \mathcal C^\infty_c(C^d_\frac34+\xi)$
such that $\tilde\f\equiv 1 $ on $C^d_{\frac12}+\xi$. By (\ref{eq:localization}) in Lemma \ref{L:localization}, 
we have that
$$
\limsup_{h\to\infty}\Big(S_{a,b}^\frac{q}{q-p}-\int\limits_{C^d_\frac34+\xi}|y|^{-d+q\rmH_{a}}|z|^{-b\frac{q}{p}}|u_h|^{q}~\!dz\Big)
\Big(\int\limits_{\R^d}|y|^{a}|z|^{-b}|\nabla(\tilde\f u_h)|^p~\!dz\Big)^{\frac{q}{q-p}}\le 0.
$$
Thanks to (\ref{eq:small}), $\tilde\f u_h\to 0$ strongly in $\mathcal D^{1,p}(\R^d;|y|^a|z|^{-b}dz)$. 
Thus also $\tilde\f u_h\to 0$ in $L^p(\R^d;|y|^a|z|^{-b-p}dz)$ by the 
Hardy inequality in (\ref{eq:CMN}). Since $\tilde\f\equiv 1$ on $C^d_\frac12+\xi$, we infer that
$$
\int\limits_{C^d_\frac12+\xi}|y|^a|z|^{-b-p}|u_h|^{p} =o(1) \qquad \text{for any $\xi\in\Sigma_0$.}
$$
A compactness argument easily gives
$$
\int\limits_{C^{d-k}_R\times C^k_\frac12}|y|^a|z|^{-b-p}|u_h|^{p}=o(1),
$$
hence (\ref{eq:included0}) implies
$$
\int\limits_{A_R}|y|^a|z|^{-b-p}|u_h|^{p}\ge \delta +o(1)~\!,
$$
where $A_R$ is the open set in (\ref{eq:Omega}).
This contradicts the Rellich Theorem and completes the proof.
\QED

\section{A closer look at the limiting cases}
\label{S:strict}

The present section is divided into two subsections. In the first one, we 
deal with the critical case $q=\pstar$ and  provide some sufficient conditions for the validity of the inequality $S_{a,b,\gamma}(\pstar)<S$. If $\gamma>b$, this guarantees 
the existence of extremals for $S_{a,b,\gamma}(\pstar)$, thanks to $ii)$ in Theorem \ref{T:Main2}. We also 
provide examples  in which the best constant 
$S_{a,b,\gamma}(\pstar)$ is not achieved.

\medskip

We emphasize that in the bottom case it could happen that $S_{a,b,b}(p^*)<S$ but $S_{a,b,b}(p^*)$ is not achieved.
This phenomenon can arise if $S_{a,b,b}(p^*)=M_a(\pstar)$, see Section \ref{S:Hilbert} for 
concrete examples.

\medskip

In the second part we concentrate our attention on the bottom case $\gamma=b$, and give two sufficient conditions  for the 
strict inequality $S_{a,b,b}(q)<M_a(q)$, which implies the existence of minimizers for $S_{a,b,b}(q)$,
see Theorem \ref{T:Main3}.

In the Hilbertian case $p=2$ a more complete description of the behavior of $b\mapsto S_{a,b,b}(q)$ will be given in Section \ref{S:Hilbert}.

\subsection{Critical case $q=\pstar$}

We  assume $p<d$ and study  the best constant $S_{a,b,\gamma}(\pstar)$ in (\ref{eq:S_critical}).
Thanks to \cite[Theorem 0.1]{GM}, we  know that in the purely spherical case we have that $S_{0,b,b}(\pstar)<S$ provided that
$b>0 $, whereas
in the purely cylindrical case it holds that $S_{a,0,0}(\pstar)<S$ if $a<0$ 
(notice that the Sobolev constant $S$ is recovered by taking $a=b=0$).

We now turn our attention to mixed weights.
For the sake of clarity, we recall the assumptions needed to have 
that $S_{a,b,\gamma}(\pstar)>0$, namely, 
$$
k+a>0~,\qquad b<p\rmH_a~,\qquad  \gamma>b-\big(k\frac{p}{\pstar}+a\big)
~,\qquad \gamma\ge b.
$$
Notice that the bottom case $\gamma=b$ is allowed if and only if $k\frac{p}{\pstar}+a>0$.

\medskip

A first set of assumptions to have $S_{a,b,\gamma}(\pstar)<S$ is obtained by observing that
$S_{a,b,\gamma}(\pstar)\to 0$ when the weight in the denominator in (\ref{eq:S_critical}) becomes not locally integrable. 

\begin{Theorem}
\label{T:simple_crit}
\begin{itemize}
\item[$i)$] If $k\frac{p}{\pstar}+a>0$, then for any $\gamma\ge p\rmH_a$, there exists $\eps=\eps(\gamma)>0$
such that $S_{a,b,\gamma}(\pstar)$ is achieved if 
$$
p\rmH_a-\eps<b<p\rmH_a;
$$

\item[$ii)$] If $k\frac{p}{\pstar}+a\le 0$ and $k+a>0$, then for any $\gamma\in \R$ there exists $\eps=\eps(\gamma)>0$ such that $S_{a,b,\gamma}(\pstar)$ is achieved if 
$$
\min\Big\{p\rmH_a, \gamma+\big(k\frac{p}{\pstar}+a\big)\Big\}-\eps<b<\min\Big\{p\rmH_a, \gamma+\big(k\frac{p}{\pstar}+a\big)\Big\}.
$$

\end{itemize}
\end{Theorem}

\proof
Let $k\frac{p}{\pstar}+a>0$ and $\gamma\ge p\rmH_a$. The weight
$|y|^{-d+\gamma \frac{\pstar}{p}}|z|^{-\gamma \frac{\pstar}{p}}$ is not locally integrable, thus  $S_{a,b,\gamma}(\pstar)\searrow 0^+$ for $b\nearrow p\rmH_a$. 
Therefore, if $b$ is close enough to $p\rmH_a$ then $S_{a,b,\gamma}(\pstar)<S$ and hence $S_{a,b,\gamma}(\pstar)$ is achieved. This proves $i)$.

\medskip
Assume now $k\frac{p}{\pstar}+a\le 0$ and $k+a>0$. Fix $\gamma\in \R$. If $\gamma\ge (d-k)\frac{p}{\pstar}$ then
$\gamma+\big(k\frac{p}{\pstar}+a\big)\ge p\rmH_a$ and one can argue as for $i)$. 

Otherwise, we have that $S_{a,b,\gamma}(\pstar)\searrow 0^+$ for $b\nearrow \gamma+\big(k\frac{p}{\pstar}+a\big)$,
since the weight
$|y|^{-k}|z|^{-\gamma \frac{\pstar}{p}}$ is not locally integrable. 
It follows that if $b$ is close enough to $\gamma+\big(k\frac{p}{\pstar}+a\big)$ then $S_{a,b,\gamma}(\pstar)<S$ and hence $S_{a,b,\gamma}(\pstar)$ is achieved. 
\QED

\medskip
The next two results involve the constant $S_{a,b,b}(\pstar)$ (bottom case $\gamma=b$). 

\begin{Theorem}
\label{T:strict_new}
Let $b<p\rmH_a$. If $a<0$ assume in addition that  
$k\frac{p}{\pstar}+a>0$.
\begin{itemize}
\item[$i)$] If $S_{a,b,b}(p^*)=S$, then for any $\gamma>b$ we have that $S_{a,b,\gamma}(\pstar)=S$ and $S_{a,b,\gamma}(\pstar)$ is not achieved.
In particular, for $a=b=0$ and for any $\gamma>0$ the infimum $S_{0,0,\gamma}(\pstar)$ is not achieved.
\item[$ii)$] If $S_{a,b,b}(p^*)<S$, then there exists $\eps>0$ such that for any $\gamma\in(b,b+\eps)$, we have  $S_{a,b,\gamma}(\pstar)<S$, and thus $S_{a,b,\gamma}(\pstar)$ is  achieved.
\end{itemize}
\end{Theorem}

\proof
Recall that $S_{a,b,\gamma}(p^*)\le S$ for any $\gamma>b$ by Theorem \ref{T:Main1}. Since the function  $\gamma\mapsto S_{a,b,\gamma}(\pstar)$ is non decreasing, we have
$$
S_{a,b,b}(p^*)\le S_{a,b,\gamma}(p^*)\le S \qquad \text{for any $\gamma> b$}.
$$ 

If $S_{a,b,b}(p^*)=S$, then  $S_{a,b,\gamma}(p^*)= S$ for any $\gamma> b$. 
The strict monotonicity of the ratio in (\ref{eq:S_critical}) for $u\neq 0$ then easily implies that $S_{a,b,\gamma}(p^*)$ is not achieved,
argue by contradiction. 

If $S_{a,b,b}(p^*)<S$ then there exists $u_0\in \mathcal C^\infty_c(\R^d)$ such that 
$$
\frac{\displaystyle\int_{\R^d}|y|^{a}|z|^{-b}|\nabla u_0|^p~\!dz}
{\Big(\displaystyle\int_{\R^d}|y|^{a\frac{\pstar}{p}}|z|^{-b\frac{\pstar}{p}}|u_0|^\pstar~\!dz\Big)^\frac{p}{\pstar}}<S.
$$
Finally, we estimate 
$$
\begin{aligned}
\lim_{\gamma\searrow b}S_{a,b,\gamma}(p^*)\le \lim_{\gamma\searrow b}&~
\frac{\displaystyle\int_{\R^d}|y|^{a}|z|^{-b}|\nabla u_0|^p~\!dz}
{\Big(\displaystyle\int_{\R^d}\big(|y||z|^{-1}\big)^{(\gamma-b)\frac{\pstar}{p}}~\!|y|^{a\frac{\pstar}{p}}|z|^{-b\frac{\pstar}{p}}|u_0|^\pstar~\!dz\Big)^\frac{p}{\pstar}}\\
&
= \frac{\displaystyle\int_{\R^d}|y|^{a}|z|^{-b}|\nabla u_0|^p~\!dz}
{\Big(\displaystyle\int_{\R^d}|y|^{a\frac{\pstar}{p}}|z|^{-b\frac{\pstar}{p}}|u_0|^\pstar~\!dz\Big)^\frac{p}{\pstar}}<S,
\end{aligned}
$$
which implies that $S_{a,b,\gamma}(p^*)<S$ if $\gamma>b$ is close enough to $b$.
The proof is complete.
\QED

\begin{Theorem}
\label{T:AN2}
Let $b<p\rmH_a$. If $a<0$ assume in addition that $k\frac{p}{\pstar}+a>0$.
\begin{itemize}
\item[$i)$] If $0\neq a\le b$ then there exists $\eps>0$ such that $S_{a,b,\gamma}(\pstar)$ is achieved for 
$\gamma\in (b, b+\eps)$. 
\item[$ii)$] If $0=a<b$ then there exists $\eps>0$ such that $S_{0,b,\gamma}(\pstar)$ is achieved for 
$\gamma\in [b, b+\eps)$. 
\end{itemize}
\end{Theorem}

\proof
Let $u\in \mathcal C^\infty_c(\R^d)\setminus\{0\}$ be radially symmetric. 
By using spherical coordinates it is easy to see that 
$$
\frac{\displaystyle\int_{\R^d}|y|^{a}|z|^{-b}|\nabla u|^p~\!dz}
{\Big(\displaystyle\int_{\R^d}|y|^{a\frac{\pstar}{p}}|z|^{-b\frac{\pstar}{p}}|u|^\pstar~\!dz\Big)^\frac{p}{\pstar}} = G_a
\frac{\displaystyle\int_{\R^d}|z|^{a-b}|\nabla u|^p~\!dz}
{\Big(\displaystyle\int_{\R^d}|z|^{(a-b)\frac{\pstar}{p}}|u|^\pstar~\!dz\Big)^\frac{p}{\pstar}}~,
$$
where 
$$
G_a=\frac{\displaystyle{\fint_{\S^{d-1}}\Ps^a d\sigma}}
{\displaystyle{\Big(~\fint_{\S^{d-1}}\Ps^{a\frac{\pstar}{p}} d\sigma\Big)^\frac{p}{\pstar}}}~\!,
$$
and $\Pi:(x,y)\mapsto y$ is the orthogonal projection   on $\R^{k}$.

Notice that $G_a<1$ by H\"older inequality, unless $a=0$. 
Next,
Lemma 3.1 in   \cite{Ho} gives 
$$
\inf_{u\in \mathcal D^{1,p}(\R^d;|z|^{a-b}dz)\atop u=u(|z|)}\frac{\displaystyle\int_{\R^d}|z|^{a-b}|\nabla u|^p~\!dz}
{\Big(\displaystyle\int_{\R^d}|z|^{(a-b)\frac{\pstar}{p}}|u|^\pstar~\!dz\Big)^\frac{p}{\pstar}}
= 
\Bigg[\frac{d-p+a-b}{d-p}\Bigg]^{p-\frac{p}{d}} S~\!.
$$
 Therefore, we obtain the estimate
\begin{equation}
\label{eq:old_lemma}
S_{a,b,b}(\pstar)\le G_a \Bigg[\frac{d-p+a-b}{d-p}\Bigg]^{p-\frac{p}{d}} S.
\end{equation}

Assume $a\neq 0$, so that $G_a<1$. Then (\ref{eq:old_lemma}) implies that 
$S_{a,b,b}(\pstar)<S$ for any $b\ge a$, thus $i)$ follows by $ii)$ in Theorem \ref{T:strict_new}.

\medskip

If $b>a=0$ we know that  $S_{0,b,b}(\pstar)$ is achieved by \cite[Theorem 0.2]{GM}.
Further, we have $G_0=1$ but (\ref{eq:old_lemma}) again yields the strict inequality $S_{0,b,b}(\pstar)<S$.
Thus $S_{0,b,\gamma}(\pstar)$ is achieved if $\gamma>b$ is close enough to $b$, use $ii)$ in Theorem \ref{T:strict_new} once more.
\QED

\subsection{Bottom case $\gamma=b$}
\label{SS:Bottom}

In the first result we get a condition to have $S_{a,b,b}(q)<M_a(q)$ by  exploiting 
the integrability properties of the weight in the denominator of (\ref{eq:Mb}).

\begin{Theorem}
\label{T:simple}
Let (\ref{eq:assu_cyl}) holds. Then there exists $\eps=\eps(a,q)>0$ such that 
$S_{a,b,b}(q)$ is achieved for any $b\in (p\rmH_a-\eps,p\rmH_a)$.
\end{Theorem}

\proof The weight
$|y|^{-d+q\rmH_a}|z|^{-q\rmH_a}$ is not locally integrable, thus  $S_{a,b,b}(q)\searrow 0^+$ for $b\nearrow p\rmH_a$. 
Therefore, if $b$ is close enough to $p\rmH_a$ then $S_{a,b,b}(q)<M_a(q)$ and hence $S_{a,b,b}(q)$ is achieved
by Theorem \ref{T:Main3}.
\QED

Due to the lack of information about the minimization problem (\ref{eq:Ma}) (purely cylindrical case),
in addition to  Theorem \ref{T:simple} we can only handle the case
$a=0$. Recall that $S_{0,b,b}(\pstar)<M_0(\pstar)=S$ is achieved for any $b\in(0,p\rmH_0)$ by $ii)$ in Theorem \ref{T:AN2}). In the subcritical case the following result holds.

\begin{Theorem}
\label{T:AN}
Let $a=0$ and assume  
$$p(d-k)<q(d-p)<dp.$$
Then there exists $\eps\in(0,p\rmH_0]$ such that $S_{0,b,b}(q)$ is achieved for any $b\in(0,\eps)$.
\end{Theorem}

\proof
Since $q<\pstar$, the infimum $M_0(q)$ is achieved by a  function $\omega\in\mathcal D^{1,p}(\R^d)$, see  \cite{GM}. 
Using symmetric rearrangements first in the $x$-variable and then in the $y$-variable (as in \cite[Theorem 3.1]{SSW}) 
we can assume that $\omega$ is positive and cylindrically symmetric. In particular,  
$\omega=\omega(|x|,|y|)$ is  decreasing in both $|x|$ and  $|y|$.
Setting $v=|z|^{b/p}\omega$, we have
$$
S_{0,b,b}(q)\le 
\frac{\displaystyle{\int_{\R^d} |z|^{-b}|\nabla v|^p~\!dz}}
{\displaystyle{\Big(\int_{\R^d} |y|^{-d+q\rmH_0}|z|^{-b\frac{q}{p}}|v|^q~\!dz\Big)^\frac{p}{q}}}
~\!=~\!\frac{\displaystyle{\int_{\R^d} \big|\nabla \omega +\frac bp \omega z|z|^{-2}\big|^p~\!dz}}
{\displaystyle{\Big(\int_{\R^d} |y|^{-d+q\rmH_0}|\omega|^q~\!dz\Big)^\frac{p}{q}}}~\!=:f(b).
$$
Notice that $f(b)<\infty$, thanks to the Hardy and Maz'ya inequalities. Further, $f(0)=M_0(q)$, $f$ is differentiable and
$$
f'(0)=c_\omega \ird |\nabla \omega|^{p-2} (\nabla \omega\cdot z) \omega |z|^{-2}~\!dz <0~\!.
$$
Thus $S_{0,b,b}(q)<M_0(q)$ for $b>0$ small enough, which concludes the proof. 
\QED

\begin{Remark} 
\label{R:symmetry}
We conjecture that the Maz'ya constant $M_a(q)$ is achieved by a radially decreasing minimizer
whenever $a\le 0$. This could allow us to extend Theorem  \ref{T:AN} 
to  parameters $a<0$. For $p=2$ this is the case, see \cite[Corollary 5.2]{GM2}.
\end{Remark}

\section{The Hilbertian, bottom case}
\label{S:Hilbert}

Here we take $p=2<q$. The hypotheses in (\ref{eq:assu_cyl}) can be written as follows:
\begin{equation}
\label{eq:2_assu}
\begin{aligned}
&a>\frac{2}{q}\qquad&&\text{if $d=2$,}\\
&a+k\frac{2}{2^*}>2(d-k)\Big(\frac{1}{q}-\frac{1}{2^*}\Big)\ge 0
\qquad &&\text{if $d\ge 3$.}
\end{aligned}
\end{equation}

We introduce the functional transform
\begin{equation}
\label{eq:Ta}
(T_bu)(z)=|z|^{-\frac{b}{2}}u(z)~,\qquad u\in \mathcal C^\infty_c(\R^d\setminus\{0\})~\!.
\end{equation}
Integration by parts easily gives
\begin{equation}
\label{eq:T}
\int\limits_{\R^d}|y|^{a}|z|^{-b}|\nabla u|^2~\!dz=
\int\limits_{\R^d}|y|^{a}|\nabla (T_bu)|^2~\!dz+\tfrac{b}{2}\big(\tfrac{b}{2}-2\rmH_a\big)\int\limits_{\R^d}|y|^{a}|z|^{-2}|T_bu|^2~\!dz,
\end{equation}
recall that $\rmH_a=\frac{d-2+a}{2}$. 
The assumptions in (\ref{eq:2_assu}) imply that $\rmH_a>0$, so that 
the following Hardy--type inequality holds:
\begin{equation}
\label{eq:Hardy2}
\int\limits_{\R^d}|y|^{a}|\nabla v|^2~\!dz\ge \rmH_a^2 \int\limits_{\R^d}|y|^{a}|z|^{-2}|v|^2~\!dz~\!\quad\text{for any $v\in \mathcal D^{1,2}(\R^d;|y|^a dz)$,}
\end{equation}
compare with (\ref{eq:CMN}). In turn, (\ref{eq:Hardy2}) implies that the inequality
\begin{equation}
\label{eq:Hardy3}
\int\limits_{\R^d}|y|^{a}|\nabla v|^2~\!dz+\tfrac{b}{2}\big(\tfrac{b}{2}-2\rmH_a\big)\int\limits_{\R^d}|y|^{a}|z|^{-2}|v|^2~\!dz
\ge \rmH_{a-b}^2 \int\limits_{\R^d}|y|^{a}|z|^{-2}|v|^2~\!dz
\end{equation}
holds for any $v\in \mathcal D^{1,2}(\R^d;|y|^a dz)$  and any $b\in\R$.

Next, let $b<2\rmH_a$, so that $\rmH_{a-b}>0$. Then 
$\mathcal C^\infty_c(\R^d\setminus\{0\})$ is dense
in $\mathcal D^{1,2}(\R^d;|y|^a|z|^{-b} dz)$ and in $\mathcal D^{1,2}(\R^d;|y|^a dz)$, use for instance \cite[Lemma 2.1]{CMN}.
Therefore, (\ref{eq:T}) and (\ref{eq:Hardy3}) imply that  the transform
$$T_b: \mathcal D^{1,2}(\R^d;|y|^a|z|^{-b} dz)\to \mathcal D^{1,2}(\R^d;|y|^a dz)$$ 
is an invertible isometry, provided that  $\mathcal D^{1,2}(\R^d;|y|^a dz)$
is endowed with the equivalent Hilbertian norm given by the left hand side of (\ref{eq:Hardy3}). Thus, we have
\begin{equation}
\label{eq:FT_bottom}
\begin{gathered}
S_{a,b,b}(q)=
\inf_{v\in \mathcal D^{1,2}(\R^d;|y|^a dz)}
\mathcal J_b(v)~\!,
\end{gathered}
\end{equation}
where
$$
\mathcal J_b(v)= \frac{\displaystyle~
\int_{\R^d}|y|^{a}|\nabla v|^2~\!dz+\tfrac{b}{2}\big(\tfrac{b}{2}-2\rmH_a\big)\int_{\R^d}|y|^{a}|z|^{-2}|v|^2~\!dz~}
{\displaystyle{\Big(
\int_{\R^d}|y|^{-d+q\rmH_{a}}|v|^q~\!dz\Big)^\frac{2}{q}}}~\!.
$$
Moreover, the minimization problems in (\ref{eq:Mb2}) and  (\ref{eq:FT_bottom}) are equivalent.

\medskip

Before going further, we recall that the Maz'ya constant $M_a(q)$ is achieved if $q(d-2)<2d$; for the critical case
$d\ge 3$, $q=2^*$ see Appendix \ref{A}.
The next result includes Theorem \ref{T:MainH}.

\begin{Theorem}
\label{T:H_bottom_sub}
Let $p=2<q$ and let (\ref{eq:2_assu}) hold. 
\begin{itemize}
\item[$i)$] There exists $b_*\in [0,2\rmH_a)$ such that 
$$
\begin{cases}
S_{a,b,b}(q)= M_{a}(q)&\text{if $b<b_*$ and in this case $S_{a,b,b}(q)$  is not achieved,}\\
S_{a,b,b}(q)< M_{a}(q)&\text{if $b_*<b<2\rmH_a$  and in this case $S_{a,b,b}(q)$  is achieved.}
\end{cases}
$$
In particular, $b_*=0$ if $M_a(q)$ is achieved;
\item[$ii)$] The map $b\mapsto S_{a,b,b}(q)$ is strictly decreasing and continuous on $[b_*,2\rmH_a)$.
\end{itemize}
\end{Theorem}

\proof
In the present proof we put
$$
S_b:=S_{a,b,b}(q),
$$ 
so that $S_0=M_a(q)$.

For fixed $v\in  \mathcal D^{1,2}(\R^d;|y|^a dz)\setminus\{0\}$ the function $b\mapsto \mathcal J_b(v)$ is strictly decreasing.
Thus  $b\mapsto S_{b}$ is non increasing. 

Since $S_{b}\le M_a(q)$ by Theorem \ref{T:Main3}, we have that $S_{b}= M_a(q)$ for $b<0$. 
The strict monotonicity of $b\mapsto \mathcal J_b(v)$ for $v$ fixed then implies that $S_{b}$ is not achieved for $b<0$,
argue by contradiction. 

To conclude the proof of $i)$ recall that $S_{b}<M_a(q)$ for $b$ close to $2\rmH_a$ by Theorem \ref{T:simple}. Thus 
$$
b_*=\inf\{b<2\rmH_a~|~ S_{b}<M_a(q)\}
$$
enjoys the properties in $i)$, thanks to Theorem \ref{T:Main3}.

\medskip

Since $S_b$ is achieved if $b\in (b_*,2\rmH_a)$, then the function 
$b\mapsto S_{b}$ is strictly decreasing on $[b_*,2\rmH_a)$, argue by contradiction.

\medskip

To conclude the proof we only need to show that the function 
$b\mapsto S_{b}$ is continuous on $[b_*,2\rmH_a)$. 
We start  
by pointing out the elementary inequalities
\begin{equation}
\label{eq:trivial}
0< \frac{b_1}{2}\Big(\frac{b_1}{2}-2\rmH_a\Big)-\frac{b_2}{2}\Big(\frac{b_2}{2}-2\rmH_a\Big)\le 2\rmH_a(b_2-b_1)~\!,
\end{equation}
which hold for $0\le b_1<b_2<2\rmH_a$.

Fix $b\in(b_*,2\rmH_a)$,   take a minimizer  
$u_{b}\in \mathcal D^{1,2}(\R^d;|y|^a|z|^{-b}dz)$ for $S_b$ and consider the function 
$v_b=T_b(u_b)\in  \mathcal D^{1,2}(\R^d;|y|^a dz)$. Notice  that $\mathcal J_b(v_b)=S_b$. Then take $\eps\in(b_*,b)$ and use (\ref{eq:trivial}) with $b_1=b-\eps, b_2=b$ to estimate
$$
S_{b-\eps}> S_{b}=\mathcal J_{b-\eps}(v_b)-\big(\mathcal J_{b-\eps}(v_b)-\mathcal J_{b}(v_b)\big)\ge \mathcal J_{b-\eps}(v_b)-  c(v_b)\eps,
$$
where $c(v_b)$ only depends on $\rmH_a$ and on the weighted $L^2$ and $L^q$ norms of $v_b$. Since $ \mathcal J_{b-\eps}(v_b)\ge S_{b-\eps}$ we infer
$$
0< S_{b-\eps}-S_{b}\le  c(v_b)~\!\eps~\!.
$$
We proved  the left continuity of $b\mapsto S_{b}$ on $(b_*,2\rmH_a)$. 

For the right continuity fix $b\in [b_*,2\rmH_a)$ and take a positive number $\eps<\rmH_a-\frac12 b$. 
Take a minimizer  
$u_{\eps}\in \mathcal D^{1,2}(\R^d;|y|^a|z|^{-b-\eps}dz)$ for $S_{b+\eps}$ and let  
$v_\eps=T_{b+\eps}(u_\eps)\in  \mathcal D^{1,2}(\R^d;|y|^a dz)$, so that 
$\mathcal J_{b+\eps}(v_\eps)=S_{b+\eps}$. We use (\ref{eq:Hardy3}) with $b+\eps$ instead of $b$
to estimate 
$$
\rmH_{a-b-\eps}^{2}\frac{\displaystyle\int_{\R^d}|y|^{a}|z|^{-2}|v_\eps|^2~\!dz}
{\Big(\displaystyle\int_{\R^d}|y|^{-d+q\rmH_{a}}|v_\eps|^q~\!dz\Big)^\frac{2}{q}}\le \mathcal J_{b+\eps}(v_\eps)=  S_{b+\eps}< S_b.
$$
Notice that $4\rmH_{a-b-\eps}>2\rmH_a-b>0$.
Using also (\ref{eq:trivial}) with $b_1=b$ and $b_2=b+\eps$ we have 
$$
\mathcal J_{b}(v_\eps)-\mathcal J_{b+\eps}(v_\eps)\le 2\rmH_{a}~\!\eps
\frac
{\displaystyle\int_{\R^d}|y|^{a}|z|^{-2}|v_\eps|^2~\!dz}
{\Big(\displaystyle\int_{\R^d}|y|^{-d+q\rmH_{a}}|v_\eps|^q~\!dz\Big)^{\frac{2}{q}}}\le 2^5~\!\frac{ S_b\rmH_a}{(2\rmH_a-b)^2} ~\!\eps~\!.
$$
Thus 
$$
\begin{aligned}
S_b> S_{b+\eps}&=\mathcal J_{b+\eps}(v_\eps)= \mathcal J_{b}(v_\eps)-\big(\mathcal J_{b}(v_\eps)-\mathcal J_{b+\eps}(v_\eps)\big)\ge S_{b}- c~\!\eps~\!,
\end{aligned}
$$
which gives the right continuity of $b\mapsto S_{b}$ on $[b_*,2\rmH_a)$
and concludes the proof.
\QED

In the subcritical case $q<2^*$ the Maz'ya constant $M_a(q)$ is achieved, hence $b_*=0$. In the critical case
$d\ge 3$, $q=2^*$ our informations on $b_*$ are not complete yet, due to the fact that there are still few open questions concerning the 
Maz'ya constant $M_a(2^*)$, see Appendix \ref{A}.

\begin{Theorem}
\label{T:almost}
Let $d\ge 3$. If $a< 0$ assume that $k\frac{2}{\twostar}+a>0$. Then
the infimum $S_{a,b,b}(\twostar)$ is achieved for any $b\in[0,2\rmH_a)$ (thus, $b_*=0$)
if one of the following conditions hold:
$$
a\le 0\qquad \text{or}\qquad \text{$d\ge 4$, $k=1$ and $a\in(0,2)$.}
$$
\end{Theorem}

\proof 
In both cases, the Maz'ya constant $M_a(2^*)$ is achieved; see \cite[Appendix A.2]{GM} (see also Appendix~\ref{A} for further details), 
thus the statement is an immediate consequence
of Theorem \ref{T:H_bottom_sub}.
\QED

Combined, Theorem~\ref{T:almost} and the next result give a clear picture of the behavior of the map $b\mapsto S_{a,b,b}(\twostar)$, leaving open only the case $d=3$, $k=1$ and $a\in(0,1)$ (recall that $S_{0,0,0}(\twostar)$ is the Sobolev constant).

\begin{Theorem}
\label{T:H_ne}
Let $d\ge 3$ and let $a\ge 0$, $b\le 0$ and assume that $a\neq 0$ or $b\neq 0$. Then for any $\gamma\ge b$ 
we have that $S_{a,b,\gamma}(\twostar)=S$ if one of the following conditions hold:

$H1)$ $a=0$;

$H2)$ $k\ge 2$ and $a>0$;

$H3)$ $d\ge 4$, $k=1$  and $a\ge 2$;

$H4)$ $d=3$, $k=1$  and $a\ge 1$.

\noindent
In all of these cases, the infimum $S_{a,b,\gamma}(2^*)$ is not achieved, for any $\gamma \ge b$.

\end{Theorem}

\proof
If $a$ satisfies one of the assumptions $H1)$ -- $H4)$ then by the available results on $M_a(\twostar)$, see Appendix \ref{A},
we have that $M_a(\twostar)=S$ and is not achieved, unless $a=0$. Therefore, 
if $b=0$ and $a> 0$, we have that $S_{a,0,0}(2^*)=M_a(2^*)$ is not achieved.

If $b<0$ we have that $S_{a,b,b}(\twostar)=S$ and $S_{a,b,b}(\twostar)$ is not achieved by $i)$ in Theorem \ref{T:H_bottom_sub}.

To conclude the proof for $\gamma>b$  use $i)$ in Theorem \ref{T:strict_new}.
\QED

\appendix 

\section{Examples}
\label{A:Examples}

Thanks to the flexibility in the choice of parameters, the class of dilation-invariant inequalities we can handle is rather broad. 
We present below a few representative examples.
It is understood that the constants on the left-hand sides are  positive.

\bigskip

\noindent
{\bf Example 1.} 
Here we take $d=2$, hence $k=1$, and $b=0$. 
For $q>2$,    $a>0$ and for any $\gamma\ge 0$ such that $a+\gamma>\frac{2}{q}$,  it holds that
$$
S_{a,0,\gamma}(q) \Big(\int\limits_{\R^2}\frac {|y|^{-2+q\frac{a+\gamma}{2}}}{(x^2+|y|^2)^{q\frac{\gamma}{4}}}~\!|u|^{q}~\!dxdy\Big)^\frac{2}{q}\le
\int\limits_{\R^2}|y|^{a}|\nabla u|^2~\!dxdy~,\quad
u\in \mathcal D^{1,2}(\R^2;|y|^adz).
$$
It is known that $S_{a,0,0}(q)=M_a(q)$ is achieved, see  Appendix \ref{A}. Thanks to 
Theorem \ref{T:Main2}, we 
also have that 
the best constant $S_{a,0,\gamma}(q)$ is achieved in $\mathcal D^{1,2}(\R^2;|y|^adz)$ 
for any $\gamma> 0$. It follows that the function $\gamma\mapsto S_{a,0,\gamma}(q)$
is strictly increasing.

\bigskip

\noindent
{\bf Example 2.}
Let $d\ge 3$ and take $a=2-k$. For future convenience we recall that $\mathcal C^\infty_c(\R^d\setminus\Sigma_0)$ is dense in 
$\mathcal D^{1,2}(\R^d;|y|^{2-k}~\!dz)$ by \cite[Lemma 2.1]{CMN}. 

Thanks to Theorem \ref{T:Main1}, for any $\gamma\ge 0$ it holds that
\begin{equation}
\label{eq:2-ka}
S_{2-k,0,\gamma}(\twostar) \Big(\ird |y|^{(2-k+\gamma)\frac{\twostar}{2}}|z|^{-\gamma\frac{\twostar}{2}}~\!|u|^{\twostar}~\!dz\Big)^\frac{2}{\twostar}\le
\ird|y|^{2-k}|\nabla u|^2~\!dz~, \quad u\in \mathcal D^{1,2}(\R^d;|y|^{2-k}~\!dz)~\!.
\end{equation}
We now distinguish three different situations, depending on the dimension of the singular set $\Sigma_0$.

\medskip
\noindent 
\begin{itemize}
\item[$-$] If $d>k\ge 3$, the existence of a minimizer for $S_{2-k,0,0}(\twostar) =M_{2-k}(\twostar)$ has been proved in \cite[Theorem 1.1]{TT}. Further,  part $i)$ in Theorem \ref{T:AN2} provides the existence of $\eps>0$ such that  
$S_{2-k,0,\gamma}(\twostar)$ is achieved for any 
$\gamma\in (0,\eps)$. It follows that the function $\gamma\mapsto S_{2-k,0,\gamma}(\twostar)$ is strictly increasing in a right neighbourhood of $0$, and
$S\ge S_{2-k,0,\gamma}(\twostar)>M_{2-k}(\twostar)$ for any $\gamma>0$.

Next, we use the 
functional change $v=|y|^\frac{2-k}{2}u$ to get
the following improvement of the Hardy inequality (see  \cite[Corollary 2.1.3]{Ma} for $\gamma=0$):
\begin{equation}
\label{eq:2-k}
S_{2-k,0,\gamma}(\twostar) \Big(\ird |y|^{\gamma\frac{\twostar}{2}}|z|^{-\gamma\frac{\twostar}{2}}~\!|v|^{\twostar}~\!dz\Big)^\frac{2}{\twostar}\le
\ird\big(|\nabla v|^2-\Big(\frac{k-2}{2}\Big)^2  |y|^{-2}|v|^2\big)~\!dz~\!.
\end{equation}
Here $v$ runs in the Hilbert space  $\mathcal D_\#^{1,2}(\R^d)$, whose norm
is given by the right hand side of (\ref{eq:2-k}). 

Since $\big(\frac{k-2}{2}\big)^2$ is the sharp constant in the 
Hardy inequality with cylindrical weights, we see that 
$\mathcal D_\#^{1,2}(\R^d)\supsetneq \mathcal D^{1,2}(\R^d)$.
By the density result mentioned earlier, the
transform $u\mapsto  |y|^\frac{2-k}{2}u$ is an invertible isometry 
$\mathcal D^{1,2}(\R^d;|y|^{2-k} dz)\to \mathcal D_\#^{1,2}(\R^d)$. As a consequence,
the best constant in (\ref{eq:2-k}) is achieved in $\mathcal D_\#^{1,2}(\R^d)$ provided that $\gamma>0$ is 
sufficiently small. For $\gamma=0$
see \cite{TT}.

\item[$-$] If $k=2$ then $S_{2-k,0,\gamma}(\twostar)=S$ and $S_{2-k,0,\gamma}(\twostar)$ is not 
 achieved if $\gamma>0$ by $i)$ in Theorem \ref{T:strict_new}.

\item[$-$] Lastly, let $k=1$, so that $\R^d\setminus \Sigma_0$ is the disjoint union of 
the half-spaces $\R^d_\pm:=\{\pm y>0\}$. Since
$\mathcal C^\infty_c(\R^d\setminus\Sigma_0)$ is dense in 
$\mathcal D^{1,2}(\R^d;|y|~\!dz)$, any function in $\mathcal D^{1,2}(\R^d;|y|~\!dz)$
is the sum of two functions both vanishing on a half-space. We see that
 (\ref{eq:2-ka}) transforms into 
$$
S_{1,0,\gamma}(\twostar) \Big(\int\limits_{\R^d_+} y^{(1+\gamma)\frac{\twostar}{2}}|z|^{-\gamma\frac{\twostar}{2}}~\!|u|^{\twostar}~\!dz\Big)^\frac{2}{\twostar}\le
\int\limits_{\R^d_+} y|\nabla u|^2~\!dz~\!,\quad u\in \mathcal D^{1,2}(\R^d_+;y~\!dz).
$$
The constant in the left hand side is sharp, use  a convexity argument as in \cite[Lemma 1.2]{Mu}.

\medskip

If $d\ge 4$, it is known that  $M_{1}(\twostar)<S$ and that  $M_{1}(\twostar)$ is achieved, see Appendix \ref{A}. Thus
$M_{1}(\twostar)<S_{1,0,\gamma}(\twostar)<S$ and $S_{1,0,\gamma}(\twostar)$ is achieved for 
$\gamma>0$ small enough by $ii)$ in Theorem \ref{T:strict_new}. 

If $d=3$, by an unexpected nonexistence result in \cite{BFL} (see also \cite{MS} and Appendix \ref{A}),
we have that $M_1(\twostar)=S$, and is not achieved. Therefore 
$S_{1,0,\gamma}(\twostar)= S$ is not achieved for any $\gamma\ge 0$, use Theorem \ref{T:strict_new}.

As before,  via the functional change $v=y^\frac12 u$, we find the equivalent inequality
with sharp constants 
\begin{equation}
\label{eq:d=d+}
S_{1,0,\gamma}(\twostar) \Big(\int\limits_{\R^d_+} y^{\gamma\frac{\twostar}{2}}|z|^{-\gamma\frac{\twostar}{2}}~\!|v|^{\twostar}~\!dz\Big)^\frac{2}{\twostar}\le
\int\limits_{\R^d_+} \big(|\nabla v|^2-\frac14  y^{-2}|v|^2\big)~\!dz~,\quad
v\in \mathcal D^{1,p}_\#(\R^d_+)~\!,
\end{equation}
where $\mathcal D_\#^{1,2}(\R^d_+)\supsetneq \mathcal D^{1,2}(\R^d_+)$ is the completion of $C^\infty_c(\R^d_+)$ with respect to the norm given by the
right hand side of (\ref{eq:d=d+}). 

\end{itemize}

\noindent
{\bf Example 3.} 
Here we focus our attention on Hardy-Sobolev type inequalities which are
related to the Hardy  inequalities investigated in \cite{CMN}.

For $a,b$ satisfying (\ref{eq:CMN_assu}) and $q>p$ we choose
\begin{equation}
\label{eq:gamma}
{\gamma}= b+\frac1q\big((d+a)p-q(d-p+a)\big)
\end{equation}
to obtain inequalities of the form
\begin{equation}
\label{eq:CMN_HS}
S_{a,b,\gamma}(q)
\Big(\int\limits_{\R^d}|y|^{a}|z|^{-{\gamma}\frac{q}{p}} |u|^q~\!dz\Big)^\frac{p}{q}\le \int\limits_{\R^d}|y|^{a}|z|^{-b}|\nabla u|^p~\!dz~\!, \quad 
u\in \mathcal D^{1,p}(\R^d;|y|^{a}|z|^{-b}~\!dz)~\!.
\end{equation}
We know that $S_{a,b,\gamma}(q)>0$ if and only if (\ref{eq:assu}) hold. 
For $\gamma$ given by  (\ref{eq:gamma}), the condition $\gamma\ge b$
translates into the constraint
$q(d-p+a)\le (d+a)p$. 
Thus we can  combine the assumptions $a1)$ and $a2)$ in (\ref{eq:assu}) into the single requirement
$$
q(d+a_+-p)\le (d+a_+)p~\!,
$$
where $a_+=\max\{a,0\}$.
In other words, in dealing with (\ref{eq:CMN_HS}), the {\it effective dimension}  is $d+a_+$; if $d+a_+>p$ then the {\it critical exponent} is
$$
p^*_{d+a_+}:=\dfrac{(d+a_+)p}{d+a_+-p}~\!.
$$
Let us briefly comment on the underlying phenomena.

When $a< 0$, $p<d$ and $q=\pstar$, we have the inequality
\begin{equation}
\label{eq:HS_CMN+}
S_{a,b,b-a\frac{p}{d}}(\pstar)
\Big(\int\limits_{\R^d}|y|^{a}|z|^{-(b-a\frac{p}{d})\frac{\pstar}{p}}|u|^\pstar~\!dz\Big)^\frac{p}{\pstar}\le \int\limits_{\R^d}|y|^{a}|z|^{-b}|\nabla u|^p~\!dz~\!, \quad 
u\in \mathcal D^{1,p}(\R^d;|y|^{a}|z|^{-b}~\!dz)
\end{equation}
The groups of translations and dilations in $\R^d$ 
may be responsible of the nonexistence of minimizers in case $S_{a,b,b-a\frac{p}{d}}(\pstar)=S$,
because minimizing sequences might concentrate at points $z\notin \Sigma_0$.
Notice that $\gamma=b-a\frac{p}{d}>b$ in this case.

When $a > 0$, $d + a > p$, and $q = p^*_{d+a}< \pstar$, 
the problem falls into the \textit{bottom case}: the only responsible 
of the lack of compactness are 
translations along the singular set $\Sigma_0$. 
The corresponding limiting inequality is
\begin{equation}
\label{eq:HS_CMN-}
S_{a,b,b}(p^*_{d+a})
\Big(\int\limits_{\R^d}|y|^{a}|z|^{-b\frac{p^*_{d+a}}{p}} |u|^{p^*_{d+a}}~\!dz\Big)^\frac{p}{p^*_{d+a}}\le \int\limits_{\R^d}|y|^{a}|z|^{-b}|\nabla u|^p~\!dz~\!, \quad 
u\in \mathcal D^{1,p}(\R^d;|y|^{a}|z|^{-b}~\!dz)
\end{equation}
Minimizers might not exist if 
$S_{a,b,b\frac{p}{d}}(p^*_{d+a})=M_a(p^*_{d+a})$.

\medskip

To conclude, we reformulate (\ref{eq:HS_CMN+}) in the "Caffarelli--Silvestre" \cite{CS} setting.
Notice that $\frac{2(n+1)}{n-1}$ is the critical Sobolev exponent for functions in 
$\mathcal D^{1,2}(\R^{n+1})$, $n\ge 2$.

\begin{Corollary}
\label{C:CS}
Let $n\ge  2$, $\frac{1}2\le s<1$. Then there exists a constant $c>0$ such that
$$
c
\Big(\int\limits_{\R^{n+1}}|y|^{1-2s}|z|^{\frac{2(1-2s)}{n-1}} |u|^{\frac{2(n+1)}{n-1}}~\!dz\Big)^\frac{n-1}{n+1}\le \int\limits_{\R^{n+1}}|y|^{1-2s}|\nabla u|^2~\!dz~\!, \quad 
u\in \mathcal D^{1,2}(\R^{n+1};|y|^{1-2s}dz)~\!.
$$
\end{Corollary} 

\begin{Remark}
With the notation in Corollary \ref{C:CS}, we have that if
$s<\frac{1}2$ (and $n\ge 1$), then  the highest exponent $q$ we can reach is $ \frac{2(n+2(1-s))}{n-2s}<\frac{2(n+1)}{n-1}$.
The corresponding Hardy-Sobolev inequality (\ref{eq:HS_CMN-}) involves purely cylindrical weights. 
\end{Remark}

\section{On the Maz'ya constants}
\label{A}

The next  result was proved in \cite{GM} in case $p<d$ (and $1\le k\le d$). The proof carries over without modification when $p\ge d$.

\begin{Proposition}
\label{P:tutto}
Assume that (\ref{eq:assu_cyl}) is satisfied. If either $q(d-p)<dp$, or  $q=\pstar$ and $M_a(\pstar) <S$, then $M_a(q)$ is achieved in  $\mathcal D^{1,p}(\R^d;|y|^adz)$.

In particular, if  $p<d$ and $a<0$, then 
$M_a(\pstar) <S$, hence $M_a(\pstar)$ is achieved.
\end{Proposition}

For general $p<d$ and $a>0$ very little is known about the Maz'ya constant $M_a(\pstar)$.
The proposition below summarizes the known results in the Hilbertian case $p=2$. For details and proofs see \cite[Proposition A.10]{GM} for items $i)$ and $ii)$, and to 
\cite[Section 6]{MS} for  $iii)$; see also \cite{BFL} for the case $a=1$.

\begin{Proposition}
\label{P:3}
Let $p=2$, $d\ge 3$ and let $a>0$. 
\begin{itemize}
\item[$i)$] If $k\ge 2$ then $M_a(\twostar)=S$ and $M_a(\twostar)$ is not achieved;
\item[$ii)$] If $k= 1$ and $d\ge 4$ then $M_a(\twostar)<S$ if and only if $0<a<2$; thus for $a\ge 2$ the Maz'ya constant $M_a(\twostar)$ is not achieved;
\item[$iii)$] If $k= 1$, $d=3$ and $a\ge 1$ then $M_a(\twostar)=S$ and the Maz'ya constant $M_a(\twostar)$ is not achieved.
\end{itemize}
\end{Proposition}

We remark that the case $d=3$, $k=1$, $0<a<1$ remains open.

\medskip

Very little is known about the symmetry properties of extremals for
$M_a(q)$ when $a\neq 0$. 
We cite \cite{CL} for a discussion on this subject, see also Remark~\ref{R:symmetry} in Subsection \ref{SS:Bottom}.

Mancini, Fabbri and Sandeep 
applied the moving plane method to show that, for $d>k\ge 2$
and $q\in(2,\twostar$), any nonnegative  solution $v\in\mathcal D^{1,2}(\R^d)$ to
$$
-\Delta v=|y|^{-d+q\rmH_0}~v^{q-1}\quad{\rm  in}~\R^d~\!,
$$
is cylindrically symmetric, see \cite[Theorem 2.1]{MFS}. This was a key step in \cite{MFS} toward the complete classification of nonnegative solutions to \eqref{eq:Ma} in the case
$q=\frac{2(d-1)}{d-2}$.
As a result, the Maz'ya constant $M_0\big(\frac{2(d-1)}{d-2}\big)$ can be explicitly computed.  

The uniqueness result in \cite{MFS} was improved in 
\cite[Theorem 1.1]{CL} for solutions $v\in  \mathcal D^{1,2}_{\rm loc}(\R^d)$, without any assumptions at infinity.
More recently it has been extended in \cite{LL} to the case $p<d$, $q=\frac{p(d-1)}{d-p}$,
under the additional assumption $k\ge 3$.

To the best of our knowledge, this and the Sobolev case $a=0$, $q=\pstar$, are the only cases where the value of $M_a(q)$ is known.

In \cite[Corollary 5.2]{GM2},  the moving plane method is used to prove that if $a\le  0$, then any 
nonnegative  solution to
$$
-\div(|y|^{a}\nabla v)=|y|^{-d+q\rmH_a}v^{q-1}\qquad \text{on $\R^d\setminus\Sigma_0$,}
\quad \ird |y|^{-d+q\rmH_a}|v|^{q}~\!dz<\infty
$$
is cylindrically symmetric. Note that $v$ may be singular on $\Sigma_0$.

Lastly, we recall that Mancini and Sandeep were the first to observe and exploit, in \cite{MS}, the  connection between variational problems on 
$\R^d$ involving cylindrical weights and certain unweighted variational problems on the hyperbolic space $\mathbb H^d$.
Nowadays, this connection is quite popular and largely used.

\bigskip

{\small
\noindent
{\bf Acknowledgement.} G. Cora is partially supported by the INdAM-GNAMPA
project n. E5324001950001. G. Cora and R. Musina
are partially supported by the Progetto PRIN 2022 - Pattern formation in nonlinear phenomena - Codice n. 20227HX33Z CUP N. G53D23001720006 ”Finanziato dall'Unione Europea, NextGeneration EU – PNRR M4 C2 I1.1".  A.I. Nazarov was supported by the Ministry of Science and Higher Education of the Russian Federation (agreement 075-15-2025-344 dated 29/04/2025 for Saint Petersburg Leonhard Euler International Mathematical Institute).

\medskip
The authors would like to thank Prof. Rupert Frank for his valuable remark, which helped us  improve the exposition.
}

\end{document}